\documentclass[aos,MSNbibl,nameyear,dvips]{arximspdf}

\doi{10.1214/14-AOS1216} %kopijuoti is PTS
\volume{42}
\issue{3}
\pubyear{2014}
\firstpage{1145}
\lastpage{1165}

\makeatletter
\def\N{\mathbb{N}}
\def\R{\mathbb{R}}
\newcommand{\widebar}{\overline}

\newtheorem{theo}{Theorem}[section]
\newproclaim{defi}{Definition}[section]
\newtheorem{prop}{Proposition}[section]
\newtheorem{lemma}{Lemma}[section]
\newtheorem{corol}{Corollary}[section]
\makeatother

\begin{document}
\begin{frontmatter}

\title{Asymptotic equivalence of nonparametric diffusion and Euler
scheme experiments}
\runtitle{Equivalence for diffusion and Euler scheme}

\begin{aug}
\author[a]{\fnms{Valentine} \snm{Genon-Catalot}\ead[label=e1]{valentine.genon-catalot@parisdescartes.fr}}
\and
\author[b]{\fnms{Catherine} \snm{Lar\'edo}\corref{}\ead[label=e2]{catherine.laredo@jouy.inra.fr}}
\runauthor{V. Genon-Catalot and C. Lar\'edo}
\affiliation{Universit\'e Paris Descartes and INRA}
\address[a]{MAP5\\
CNRS UMR 8145\\
Universit\'e Paris Descartes\\
PRES Sorbonne Paris Cit\'e\\
45, rue des Saints-P\`eres\\
75006 Paris\\
France\\
\printead{e1}} %adresu isvedimo komanda gale!
\address[b]{MIA and LPMA\\
CNRS-UMR 7599\\
INRA and Universit\'e Paris Diderot\\
INRA, 78350, Jouy-en-Josas\\
France\\
\printead{e2}}
\end{aug}

% HISTORY:
\received{\smonth{6} \syear{2013}}
\revised{\smonth{3} \syear{2014}}

% ABSTRACT
%
\begin{abstract}
We prove a global asymptotic equivalence of experiments in the sense of
Le Cam's
theory. The experiments are a continuously observed diffusion with
nonparametric drift and its Euler scheme.
% which is an autoregression model.
We focus on diffusions with nonconstant-known diffusion coefficient.
The asymptotic equivalence is
proved by constructing explicit equivalence mappings
based on random time changes. The equivalence of the discretized
observation of the diffusion and
the corresponding Euler scheme experiment is then derived. The impact
of these equivalence results is that it justifies the use
of the Euler scheme instead of the discretized diffusion process for
inference purposes.
\end{abstract}

% KEYWORDS
\begin{keyword}[class=AMS]
\kwd[Primary ]{62B15}
\kwd{62G20}
\kwd[; secondary ]{62M99}
\kwd{60J60}
\end{keyword}
\begin{keyword}
\kwd{Diffusion process}
\kwd{discrete observations}
\kwd{Euler scheme}
\kwd{nonparametric experiments}
\kwd{deficiency distance}
\kwd{Le Cam equivalence}
\end{keyword}

\end{frontmatter}

%s1 #&#
\section{Introduction}\label{sec1}

Proving global asymptotic equivalence of statistical experiments
by means of the Le Cam theory of deficiency [\citet{LCY}] is an
important issue
for nonparametric estimation problems. The interest
is to obtain asymptotic results for some experiment by means of an
equivalent one. Concretely, in the case of bounded loss functions, a
solution to a nonparametric problem in an experiment yields a
corresponding solution in an asymptotically equivalent experiment.
For instance, when minimax rates of convergence in a nonparametric
estimation problem are obtained in one experiment, the same rates hold
in a globally asymptotically equivalent experiment.
The theory also allows to prove asymptotic sufficiency of the
restriction of an experiment to a smaller \mbox{$\sigma$-}field. When
explicit transformations from one experiment to another one are
obtained, statistical procedures can be carried over from one
experiment to the other one. There is an abundant literature devoted to
establishing asymptotic equivalence results. Before considering
diffusion experiments, we recall the main contributions in this domain.
The first results concern the asymptotic equivalence of density
estimation and white noise model [\citet{NU}] and nonparametric
regression and white noise [\citet{BRL}]. These results were extended
to the equivalence of nonparametric regression with random design and
white noise [\citet{BCLZ}]. The equivalence between the observation of
$n$ independent random variables $X_i, i =1, \ldots,n$ with densities
$p(x,\theta_i)$, such that $\theta_i=f(i/n)$ and a nonparametric
Gaussian shift experiment with drift linked with $f$ is proved in \citeauthor{GRA} (\citeyear{GRA,GRA2}). In \citet{BRC}, the equivalences concern Poisson processes
with nonparametric intensity and white noise. \citet{C} considers the
equivalence of a fixed design regression in two dimensions and a
Brownian sheet process with drift. This result is extended to
regression experiments with arbitrary dimension in \citet{REIS}. The
regression model with nonregular errors yields different results, the
equivalence being with independent point Poisson processes [\citet
{MR}]. A step forward in another direction concerns the equivalence of
nonparametric autoregression and nonparametric regression [\citet
{GRANE}]. Negative results are also important such as the
nonequivalence of nonparametric regression and density or white noise
when the regression function has smoothness index $1/2$ [\citet{BZ}].
To our knowledge, the only paper studying the equivalence problem for
regression with unknown variances is \citet{CA}. More recently, the
class of studied models has been enlarged to stationary Gaussian
processes with unknown spectral density which are equivalent to white
noise [\citet{GO}]. Another direction concerns inverse problems in
regression and white noise [\citet{ME}]. Unusual rates formerly
obtained by \citet{GJ} find their mathematical understanding with the
equivalence result of \citet{REIS2} where the discretization of a
continuous Gaussian martingale observed with noise on a fixed time
interval is equivalent to a Gaussian white noise experiment with the
same unusual rate ($n^{-1/4}$ instead of $n^{-1/2}$ in the noise intensity).

Diffusion models defined by stochastic differential equations have also
been investigated. References concern nonparametric drift estimation
with known constant diffusion coefficient. \citet{GCLN} \mbox{studied} the
equivalence of a transient diffusion having positive drift and small
constant diffusion coefficient with a white noise model and other
related experiments. In the case of recurrent diffusion models, global
equivalence with Gaussian white noise no longer holds [\citet{DEH} for
null recurrent diffusions, \citeauthor{DR} (\citeyear{DR,DR2}) for ergodic scalar and
multidimensional diffusions].
ARCH-GARCH models exhibit nonstandard equivalence results when compared
to their limiting diffusion experiments. In a parametric context, \citet{W} proves the nonequivalence of the GARCH-experiment with its limiting
stochastic volatility model for the natural sampling frequencies. To
get the equivalence, suitable frequencies of observations are required
[\citet{BWZ}].

Inference for continuously observed diffusion processes is well
developed [e.g., \citet{KUT}]. As the diffusion coefficient is
identified from a continuous time observation, it is assumed to be
known and inference concerns the drift coefficient. On the contrary,
inference for discretely observed diffusions is more difficult as the
transition densities are generally untractable. Statistical procedures
based on the Euler scheme corresponding to the one-step discretization
of the diffusion have been successfully carried over to the discretized
diffusion observations. In parametric inference, we may quote \citet{Genon}, \citet{La} for small diffusion coefficient, \citet{Kess} for
positive recurrent diffusions and for nonparametric inference, \citet{H}, \citet{CGCR}. Therefore, a natural issue for understanding these
results is to prove the equivalence of the discretized observation of a
diffusion and the corresponding Euler scheme experiment. Such a result
has been proved by \citet{MN} for diffusions with small-known constant
diffusion coefficient and by \citet{DR} for positive recurrent
diffusions with constant diffusion coefficient. Our aim here is to
extend this result to the case of a nonconstant-known diffusion
coefficient using random time changes which yield models with diffusion
coefficient equal to $1$. This provides a canonical way for solving the
equivalence problem. The time changed experiment coming from the Euler
scheme does not lead to an autonomous diffusion but to an It\^o process
with predictable drift which induces the main difficulties.

More precisely, we consider the diffusion process $(\xi_t)$ given by
%
%e1 #&#
\begin{equation}
\label{eq1} d \xi_t= b(\xi_t) \,dt +\sigma(
\xi_t) \,dW_t, \qquad\xi_0= \eta,
\end{equation}
where $(W_t)_{t \geq0}$ is a Brownian motion defined on a filtered
probability space $(\Omega, {\mathcal A}, ({\mathcal A}_t)_{t \geq0},
\mathbb{P})$, $\eta$
is a real valued random variable, ${\mathcal A}_0$-measurable, $b(\cdot),
\sigma(\cdot)$ are real-valued functions defined on ${\mathbb R}$. The
diffusion coefficient $\sigma(\cdot)$ is a known nonconstant function. The
drift function $b(\cdot)$ is unknown and belongs to a nonparametric class.
The sample path of $(\xi_t)$ is continuously observed on a time
interval $[0,T]$. We also consider the discrete observation of $(\xi
_t)$ at the times $t_i= ih, i\le n$ with $n=[T/h]$. For simplicity, we
assume in what follows that $T/h$ is an integer. The Euler scheme
corresponding to~(\ref{eq1}), with sampling
interval $h$ is
%
%e2 #&#
\begin{equation}
\label{eulerscheme} Z_0 = \eta, \qquad Z_i = Z_{i-1}
+ h b(Z_{i-1}) + \sqrt{h} \sigma (Z_{i-1})\varepsilon_i,
\end{equation}
where, for $i\ge1$, $t_i= ih$ and $\varepsilon_i=
(W_{t_i}-W_{t_{i-1}})/\sqrt{h}$. For performing the comparisons, we
consider $(Z_0, Z_1, \ldots, Z_n)$ with $n=T/h$. We prove the
asymptotic equivalences assuming that $n$ tends to infinity with
$h=h_n$ and $nh_n^2=T^2/n$ tending to $0$. This includes both cases
$T=nh_n$ bounded and $T\rightarrow+\infty$. Note that, for inference
in diffusion models from discrete observations, the constraint
$nh_n^2\rightarrow0$ is the standard condition for Lipschitz drift
functions [e.g., \citet{Kess}, \citet{DR}, \citet{CGCR}]. We can
also observe that statistical procedures for estimating the drift
generally do not use the knowledge of the diffusion coefficient which
appears as a nuisance parameter. \citet{CA} did a noteworthy
improvement in this direction: he proves the asymptotic equivalence of
the regression experiment with unknown variances with an experiment
having two components, the first containing information about the
variance, the second containing information on the mean. An important
open problem which has never been tackled concerns the similar result
for diffusion processes with unknown diffusion coefficient $\sigma(\cdot)$.

%%%%%%%%%%%%%%%%%%%%%%%%%%%%%%%%%%%%%%%%%%%%%%%%%%%%%%%%%%%%%%%%%%%%%
%%%%%%%%%%%%%%%%%%%%%%%%%%%%%%%%%%%%%%%%%%%%%%%%%%%%%%%%%%%%%%%%%%%%%%%%%%%
%%%%%%%%%%%%%%%%%%%%%%%%%%%%%%%%%%%%%%%%%%%%%%%%%%%%%%%%%%%%%%%%%%%%%%%%%%%%%%%
%MODIFIE

The paper is organized as follows. Assumptions and main results are given
in Section~\ref{sec2}. Theorem~\ref{thmainresult} states the equivalence
result of~(\ref{eq1}) and~(\ref{eulerscheme}) and Corollary~\ref
{corolresultsuf} states the equivalence of the discrete observation
of the diffusion and its Euler scheme. The proof of Theorem~\ref
{thmainresult} is developed in Section~\ref{sec3}. We consider random time
changes on the diffusion and on the
Euler scheme leading to processes with diffusion coefficient equal to
$1$. %In Section 3.1,
First, the classical random time change on the diffusion
which leads to an autonomous diffusion process with drift $f=b/\sigma
^2$ and diffusion coefficient equal to $1$ is recalled (Proposition
\ref{changtemps}). We prove the exact equivalence between the
diffusion experiment~(\ref{eq1}) and the random time changed
experiment (Proposition~\ref{propdeltaeetilde}).
For the Euler scheme%(Section 3.2)
, we build a continuous time accompanying experiment
(Proposition~\ref{lemeulerdisccont}). Then we introduce a random time
change leading to a process with unit diffusion coefficient. This
process characterized in Proposition~\ref{qtildeb} has a predictable
path-dependent drift term. The exact equivalence between the
corresponding experiment and the Euler scheme experiment is proved in
Theorem~\ref{propdeltaggtilde}.
Finally, for
$n \rightarrow\infty$, the asymptotic equivalence
of the two randomly stopped\vspace*{1pt} experiments is proved (Proposition~\ref
{proprandom}) under the condition $h=h_n \rightarrow0 $, $nh_n^2
\rightarrow0$, thus completing the proof of Theorem~\ref{thmainresult}.
Concluding remarks and extensions are given in Section~\ref{sec4}. Proofs are
gathered in Section~\ref{sec5}.
%In Section~\ref{secappen}, the definition and some properties of the
%Le Cam deficiency distance $\Delta$ between statistical experiments
%are recalled,
\hyperref[secappen]{Appendix} contains a short recap on the Le Cam
deficiency distance $\Delta$ and some useful auxiliary results.

%%%%%%%%%%%%%%%%%%%%%%%%%%%%%%%%%%%%%%%%%%%%%%%%%%%%%%%%%%%%%%%%%%%%%%%%%%%%%%

%s2 #&#
\section{Assumptions and main results}\label{sec2} \label{secassumptions}
We assume that the diffusion coefficient $\sigma(\cdot)$ of~(\ref{eq1})
is known, belongs to $C^{2}({\mathbb R})$ and satisfies:
\begin{longlist}[(C)]
\item[(C)] $\forall x \in{\mathbb R},
0< \sigma_{0}^{2} \leq\sigma^{2}(x) \leq\sigma_{1}^{2},
|\sigma'(x)|+|\sigma''(x)|\le K_{\sigma}$.
\end{longlist}

The function $b(\cdot)$ is unknown and such that, for $K$ a positive
constant:

\begin{longlist}[(H1)]
\item[(H1)] $ b(\cdot) \in{\mathcal F}_K= \{b(\cdot) \in C^{1}({\mathbb
R})$ and for all $ x \in{\mathbb R},
|b(x)|+|b'(x)|\le K\} $.
\end{longlist}

The constant $K$ has to exist but may be unknown.

Condition (C) and assumption~(H1) ensure that the stochastic
differential equation~(\ref{eq1}) has a unique strong solution process
$(\xi_t)_{t \geq0}$. The assumptions on $b, \sigma$ are rather
strong but allow to shorten technical proofs. Note that (H1) and (C)
include models with or without ergodicity properties. The distribution
of the initial variable $\eta$ of~(\ref{eq1}) may be known or
unknown.

Let $C({\mathbb R}^{+},{\mathbb R})$ be the space of continuous real
functions defined on $\mathbb{R}^{+}$, and denote by $(X_t, t \geq
0)$ the canonical process of $C(\mathbb{R}^{+},\mathbb{R})$ given by
$(X_t(x)=x(t), t\ge0)$ for $x \in C(\mathbb{R}^{+},\mathbb{R})$,
${\mathcal C}_{t}^{0,X}= \sigma(X_s, s \le t)$, ${\mathcal C}_{t}^X
=\bigcap_{s>t} {\mathcal C}_{s}^{0,X}$
and ${\mathcal C}^X= \sigma( {\mathcal C}_{t}^X, t \geq0)$.
Denote by $P_b$ the distribution of $(\xi_t, t \geq0)$ defined by
(\ref{eq1}) on $( C(\mathbb{R}^{+},\mathbb{R}), {\mathcal C}^X)$ and
consider the experiment associated with the continuous observation of
the diffusion
\[
{\mathcal E}_{0} = \bigl( C \bigl({\mathbb R}^+, \mathbb{R}\bigr), {
\mathcal C}^X, (P_{b}, b \in{\mathcal F}_K)
\bigr). %
\]
If $T$ is fixed or is a $({\mathcal C}_{t}^X )$-stopping time,\vspace*{-1pt} we define
the restriction $P_{b}/_{{\mathcal C}_{T}^X}$ of $P_b$ to the $\sigma
$-field $ {\mathcal C}_{T}^X$.
The experiment associated with the continuous observation of $(\xi_t)$
stopped at $T$ is
%
%e3 #&#
\begin{equation}
\label{e0t} {\mathcal E}_{0}^T = \bigl( C \bigl({
\mathbb R}^+, \mathbb{R}\bigr), {\mathcal C}_T^X,
(P_{b}/_{{\mathcal C}_{T}^X}, b \in{\mathcal F}_K) \bigr).
\end{equation}
Consider now the Euler scheme corresponding to~(\ref{eq1}), with sampling
interval $h$, defined in~(\ref{eulerscheme}). This experiment is an
autoregression model but we have rather call it Euler scheme as it is
associated with the one-step discretization of~(\ref{eq1}). Let $(\pi
_i)_{i\geq0}$ denote the canonical projections of
$\mathbb{R}^{\mathbb{N}} \rightarrow\mathbb{R}$ given by $(\pi
_i(x)=x_i,i\ge0)$ for $x \in\mathbb{R}^{\mathbb{N}}$ and set
${\mathcal G}_n =\sigma(\pi_0,\pi_1,\ldots, \pi_n)$, ${\mathcal
G}=\sigma({\mathcal G}_n, n\ge0)$.
We\vspace*{1pt} denote by $Q_b^{h}$ the distribution of $(Z_i, i\ge0)$ defined by
(\ref{eulerscheme}) on $({\mathbb R}^{\mathbb{N}},
{\mathcal B}({\mathbb R}^{\mathbb N }))$.
For\vspace*{1pt} $N$ a \mbox{$({\mathcal G}_n )$-}stopping time, we consider the
restriction $Q_{b}^{h}/_{{\mathcal G}_{N}}$ of $Q^{h}_{b}$ to
${\mathcal G}_N $.
The experiment associated with the discrete Euler scheme $(Z_i)$ with
sampling interval $h$ stopped at $N$ is
%
%e4 #&#
\begin{equation}
\label{g} {\mathcal G}^{h, N} = \bigl(\R^{\N}, {\mathcal
G}_N, \bigl(Q_{b}^{h}/_{{\mathcal G}_{N}}, b \in{
\mathcal F}_K\bigr) \bigr).
\end{equation}
We now state the main result.
%
%th2.1 #&#
\begin{theo}\label{thmainresult}
Assume \textup{(H1)--(C)}. For deterministic $N=n$, $h=h_n$, the sequences of
experiments $({\mathcal E}_0^{nh_n})$ and $({\mathcal G}^{h_n, n} )$
are asymptotically equivalent for the Le Cam distance $\Delta$ as $n
\rightarrow\infty$, if $h_n \rightarrow0   $ and
$nh_n^2\rightarrow0$:
$
\Delta({\mathcal E}_0^{nh_n},\break {\mathcal G}^{h_n, n}) \rightarrow0$.
\end{theo}
%
%Note that in the above result, $T=nh_n$ may be either fixed or tending
%to infinity.
%We may now draw as
An important consequence is the comparison of the experiment associated
with the discrete observation $(\xi_{ih}, i\le n)$ of the diffusion
with sampling interval $h$ and the experiment ${\mathcal G}^{h, n}$.
Let $P^{h}_{b}$ denote the distribution of
$( \xi_{ih})_{i \geq0}$ defined by equation~(\ref{eq1}) on $(\mathbb
{R}^{\mathbb{N}}, {\mathcal B}(\mathbb{R}^{\mathbb{N}}))$.
For $N$ a $({\mathcal G}_n $)-stopping time, let $P_{b}^{h}/_{{\mathcal
G}_{N}}$ be the restriction of $P^{h}_{b}$ to ${\mathcal G}_N $.
The experiment associated with the discrete observations $(\xi_{ih})$
with sampling $h$ stopped at $N$ is
\[
{\mathcal E}^{h, N} = \bigl(\R^{\N}, {\mathcal
G}_N, \bigl(P_{b}^{h}/_{{\mathcal G}_{N}}, b \in{
\mathcal F}_K\bigr) \bigr). %
\]
%
%co2.1 #&#
\begin{corol}\label{corolresultsuf}
Assume \textup{(H1)--(C)}. For deterministic $N=n$, $h=h_n$, the sequences of
experiments $({\mathcal E}^{h_n, n})$ and $({\mathcal G}^{h_n,n} )$ are
asymptotically equivalent for the Le Cam distance $\Delta$ as $n
\rightarrow\infty$, if $h_n \rightarrow0$ and $nh_n^2\rightarrow0$:
$\Delta({\mathcal E}^{h_n, n},\break {\mathcal G}^{h_n, n}) \rightarrow0$.
\end{corol}
\citet{MN}, \citet{DR} proved that when $\sigma(\cdot)$ is constant and
$nh_n^2$ tends to $0$, the discrete observation $(\xi_{ih_n},\break  i\le n)$
is an asymptotically sufficient statistic for $(\xi_t, t\le nh_n)$,
that is, $\Delta({\mathcal E}_0^{nh_n},\break {\mathcal E}^{h_n,
n})\rightarrow0$. For nonconstant diffusion coefficient, the latter
asymptotic sufficiency result can be deduced using the change of function
$F(x)= \int_0^x \,du/\sigma(u)$. Therefore, applying Theorem~\ref
{thmainresult} yields the corollary.

%%%%%%%%%%%%%%%%%%%%%%%%%%%%%%%%%%%%%%%%%%%%%%%%%%%%%%%%%%%%%%%%%%%%%%%%%%%%%%%%%%%%
%s3 #&#
\section{Random time changed experiments}\label{sec3}
To deal with the nonconstant diffusion coefficient $\sigma(\cdot)$, we
define experiments obtained by random time changes. For this, set
%
%e5 #&#
\begin{equation}
\label{f} f(x)= \frac{b(x)}{\sigma^{2} (x)},\qquad L=\frac{K}{\sigma
_0^2} \biggl(1+2
\frac{K_{\sigma}\sigma_1}{\sigma_0^2}\biggr).
\end{equation}
Under \textup{(H1)--(C)}, $f$ is bounded and globally Lipschitz with constant $L$.

%s3.1 #&#
\subsection{Time change on the diffusion}\label{sec3.1} Define for $x \in C({\mathbb
R}^+,\R)$, $ t, u\ge0$,
%
%e6 #&#
\begin{equation}
\label{tau} \rho_t(x)=\int_0^{t}
\sigma^2\bigl(x(s)\bigr)\,ds,\qquad \tau_u(x)=\inf\bigl\{ t\ge0,
\rho_t(x) \ge u\bigr\}.
\end{equation}
Since $\sigma(\cdot)$ is known, the functions $ \rho_t$ and $\tau_u$ are
known as well. Therefore, one is allowed to use these functions in the
construction of Markov kernels.
By \textup{(C)}, $ \rho_{+\infty}(x) =+\infty$, $\frac{u}{\sigma_1^2}
\le\tau_u(x)\le\frac{u}{\sigma_0^2}$, $\rho_{\tau_u(x)}(x)=u$,
$\tau_{\rho_t(x)}(x)=t $. Note that $\tau_u(X) $ is a stopping
time with respect to the canonical filtration $({\mathcal C}_s^X, s\ge
0)$. We introduce now a classical time changed process.
%
%pr3.1 #&#
\begin{prop}\label{changtemps} Assume \textup{(H1)--(C)}. Let $\xi$ be the
solution of~(\ref{eq1}) and set $(\zeta_u= \xi_{\tau_u(\xi)}, u\ge
0)$ and $({\mathcal G}_u= {\mathcal A}_{\tau_u(\xi)}, u\ge0)$. Then
%
%e7 #&#
\begin{equation}
\label{defY} d\zeta_u= f(\zeta_u) \,du+ dB_u,\qquad
\zeta_0= \eta,
\end{equation}
with $(B_u)$ Brownian motion w.r.t. $({\mathcal G}_u)$
%=o the filtration $({\mathcal G}_u= {\mathcal A}_{\tau_u(\xi)}, u\ge
%0)$
which satisfies the usual conditions.
\end{prop}
The proof relies on classical tools [e.g., \citet{KA}, Chapter 3, Section 4 and Chapter 5, Section 5]
and\vspace*{2pt} implies $\delta({\mathcal E}_{0}^{\tau_a(X)},{\widetilde {\mathcal E}}_{0}^{a} ) =0$ (see \hyperref[secappen]{Appendix}). The main
difficulty lies in studying the other deficiency.
Denote by ${\widetilde P}_{b}$ the distribution of $(\zeta_u, u\ge0)$ on
$C(\R^+,\R)$. We associate to the time changed process $(\zeta_u,
u\ge0)$ an experiment with sample space $C(\R^+,\R)$. For sake of
clarity, we use a distinct notation for the canonical process and
filtration. Let $(Y_u, u\ge0)$ be defined by $Y_u(y(\cdot))=y(u)$ with
$y(\cdot) \in C(\R^+,\R)$, $({\mathcal C}_u^Y, u\ge0)$ be the associated
right-continuous canonical filtration and ${\mathcal C}^Y= \sigma
({\mathcal C}_u^Y, u\ge0)$. Set
\[
{\widetilde{\mathcal E}}_{0}= \bigl(C\bigl(\R^+,\R\bigr), {\mathcal
C}^Y, ({\widetilde P}_{b}, b \in{\mathcal
F}_K) \bigr).
\]
For $A>0$ a $({\mathcal C}_u^Y)$-stopping time, define
the experiment
\[
{\widetilde{\mathcal E}}_{0}^{A}= \bigl(C\bigl(\R^+,\R
\bigr), {\mathcal C}_A^Y, ({\widetilde
P}_{b}/_{{\mathcal C}_A^Y}, b \in{\mathcal F}_K) \bigr).
\]
Define for $y \in C(\R^+,\R)$,
%
%e8 #&#
\begin{equation}
\label{tuat} \qquad T_u(y)=\int_0^u
\frac{dv}{\sigma^2(y(v))},\qquad A_t(y)= \inf\bigl\{ u\ge 0, T_u(y) \ge
t\bigr\} =T_.(y)^{-1}(t).
\end{equation}
Thus, for all $t\ge0$, $A_t(Y)$ is a $( {\mathcal C}_u^Y)$-stopping time.
%
%pr3.2 #&#
\begin{prop}\label{propdeltaeetilde} Assume \textup{(H1)--(C)}. If $x=(x(t),
t\ge0)$, $ (y(u)=\break  x(\tau_u(x)), u\ge0)$, then $A_{t}(y)= \rho
_t(x)$, $T_u(y)= \tau_u(x)$. For $a,T$ deterministic $  \Delta
({\mathcal E}_{0}^{T},{\widetilde{\mathcal E}}_{0}^{A_T(Y)})=0$ and
$\Delta({\mathcal E}_{0}^{\tau_a(X)},{\widetilde{\mathcal E}}_{0}^{a})=0$.
\end{prop}

The experiments ${\mathcal E}_{0}$ and ${\widetilde{\mathcal E}}_{0}$ are
linked by the mapping $(x(t), t \ge0) \rightarrow(y(u)= x(\tau
_u(x)), u\ge0)$. For the stopped experiments,
noting that
$\{u, \tau_u(x)\leq T\} =\{u, u \leq A_T(y)\}$, the previous mapping
links ${\mathcal E}_{0}^{T}$ and ${\widetilde{\mathcal E}}_{0}^{A_T(Y)}$.
Similarly, the experiments ${\widetilde{\mathcal E}}_{0}$ and ${\mathcal
E}_{0}$ are linked by the mapping
$(y(u), u\ge0) \rightarrow(x(t)= y(A_t(y)), t\ge0)$ and, for stopped
experiments, noting that $\{ t, A_t(y) \leq a\}= \{t, t\leq\tau
_{a}(x)\} $, this mapping links ${\widetilde{\mathcal E}}_{0}^{a}$ and
${\mathcal E}_{0}^{\tau_a(X)}$.
%%%%%%%%%%%%%%%%%%%%%%%%%%%%%%%%%%%%%%%%%%%%%%%%%%%%%%%%%%%%%%%%%%%%%

%s3.2 #&#
\subsection{Time change on the Euler scheme}\label{sec3.2}
As the discrete Euler scheme experiment~(\ref{g}) has not the same
sample space as the diffusion experiment~(\ref{e0t}), an \mbox{essential}
tool is to use the accompanying experiment of~(\ref{g}) which is the
continuous-time Euler scheme. Given a path $x=x(\cdot) \in{\mathcal
C}({\mathbb R}^+,{\mathbb R})$
and a sampling scheme $t_i=ih, i\ge1 $, we define
the diffusion-type process $\bar{\xi}_t$,
%
%e9 #&#
\begin{equation}\label{defXbar}
d \bar{\xi}_t= \bar{b}_{h}(t, \bar{\xi})
\,dt+ \bar{\sigma}_{h}(t,\bar{\xi}) \,dW_t, \qquad\bar{\xi
}_0=\eta,
\end{equation}
with
\[
\bar{b}_{h}(t,x)=\sum_{i\geq1} b
\bigl(x(t_{i-1})\bigr){1}_{(t_{i-1},t_i]
}(t),\qquad \bar{\sigma}_{h}(t,
x)=\sum_{i\geq1} \sigma \bigl(x(t_{i-1})
\bigr){1}_{(t_{i-1},t_i] }(t).
\]
Let $ Q_b$ denote the distribution of $({\bar\xi}_t, t \geq0)$
on $ (C(\mathbb{R}^{+},\mathbb{R}), {\mathcal C}^X)$ and, for $T$ a
\mbox{$({\mathcal C}_{t}^X )$-}stopping time, $ Q_b/_{{\mathcal C}_{T}^X}$ the
restriction of $ Q_b$ to $ {\mathcal C}_{T}^X$. Set
%
%e10 #&#
\begin{equation}
\label{gbar0} {\mathcal G}_{0}^T = \bigl( C \bigl({
\mathbb R}^+, \mathbb{R}\bigr), {\mathcal C}_T^X, (
Q_b/_{{\mathcal C}_{T}^X}, b \in{\mathcal F}_K) \bigr).
\end{equation}
%
%pr3.3 #&#
\begin{prop} \label{lemeulerdisccont}
For $h>0$, $N$ a $({\mathcal G}_n)$-stopping time, the Le Cam distance
between ${\mathcal G}^{h,N}$ and
$ {\mathcal G}_{0}^{N h}$ [(\ref{g}),~(\ref{gbar0})] is equal to 0:
$\Delta({\mathcal G}^{h, N},
{\mathcal G}_{0}^{Nh }) = 0$.
\end{prop}
%
%In the case of a constant diffusion coefficient, this result is proved
%in \citet{MN}. Note that the proof uses the assumption that $\sigma(\cdot)$
%is known. \\
Let us define a time changed process associated with the continuous
Euler scheme $({\bar\xi}_t)$. The study of this time changed
process is more difficult because the drift term and the diffusion
coefficient of the continuous-time Euler scheme are time and path
dependent. Let
\[
{\bar\rho}_t(x)= \int_0^{t} {\bar
\sigma}_{h}^{2}(s, x) \,ds, \qquad \bar\tau_u(x)=
\inf\bigl\{ t\ge0, {\bar\rho }_t(x) \ge u \bigr\}.
\]
Analogously,\vspace*{1pt} ${\bar\tau}_u(X)$ is a stopping time of the
canonical filtration $ {\mathcal C}^X$. With the convention $\sum_{j=0}^{i-1}=0$ for $i=0$, we have, for $i\ge0$ and $t_i <t \le t_{i+1}$,
\[
{\bar\rho}_t(x)= {\bar\rho}_{t_i}(x) +
(t-t_i) {\bar\sigma}_{h}^2 (t_i,
x)= h \sum_{j=0}^{i-1}\sigma^2
\bigl(x(t_j)\bigr) +(t-t_i)\sigma^2
\bigl(x(t_i)\bigr).
\]
Hence, $({\bar\rho}_{t}(x),t\ge0)$ is continuous, increasing on
$\R^+$ and
maps $(t_i, t_{i+1}]$ on
$({\bar\rho}_{t_i}(x), {\bar\rho}_{t_{i+1}}(x)]$. By \textup{(C)},
$
{\bar\rho}_{+\infty}(x)=+\infty$, $u/\sigma_1^2 \le
{\bar\tau}_u(x)\le(u/\sigma_0^2)+\Delta$,
and $\{t \rightarrow{\bar\rho}_t(x)\}$, $\{u \rightarrow
{\bar\tau}_u(x)\}$ are inverse. In particular, for all $i, x$,
$t_i={\bar\tau}_{{\bar\rho}_{t_i}(x)}(x)$.
For ${\bar\xi}$ solution of~(\ref{defXbar}), set $({\widebar
{\mathcal G}}_u =
{\mathcal A}_{{\bar\tau}_u({\bar\xi})})$, and define the process
%
%e11 #&#
\begin{equation}
\label{Ybar} ({\bar\zeta}_u = {\bar\xi}_{{\bar\tau
}_u({\bar\xi})},u \ge0),
\end{equation}
which is adapted to the filtration $({\widebar{\mathcal G}}_u )$
% {\mathcal A}_{{\bar\tau}_u({\bar\xi})})$. This filtration
which satisfies the usual conditions. Denote by ${\widetilde Q}_{b}$
the distribution of $({\bar\zeta}_u)$.
%%%%%%%%%%%%%%%%%%%%%%%%%%%%%%%%%%%%%%%%%%%%%%%%%%%%%%%%
%
%pr3.4 #&#
\begin{prop} \label{qtildeb}
The process $({\bar\zeta}_u)$ defined in~(\ref{Ybar}) has unit
diffusion coefficient and drift term given by [see~(\ref{f})]:
%
%e12 #&#
\begin{equation}
{\bar f}(v)= \sum_{i\ge0} f({\bar
\zeta}_{{\bar\rho}_{t_i}
({\bar\xi})}) 1_{({\bar\rho}_{t_i}({\bar\xi}),
{\bar\rho}_{t_{i+1}}({\bar\xi})]}(v),
\end{equation}
where $({\bar\rho}_{t_i}({\bar\xi}))$ are $({\widebar
{\mathcal G}}_u)$-stopping times and so,
%w.r.t. the filtration $({\widebar{\mathcal G}}_u)$
${\bar f}(v)$ is predictable w.r.t. $({\widebar{\mathcal G}}_u)$.
\end{prop}
%
%%%%%%%%%%%%%%%%%%%%%%%%%%%%%%%%%%%%%%%%%%%%%%%%%%%%%%%%
%%%%%%%%%%%%%%%%%%%%%%%%%%%%%%%%%%%%%%%%%%%%%%%%%%%%%%%%%%%
We associate to the time changed process $({\bar\zeta}_u, u\ge
0)$ an experiment with sample space $C(\R^+,\R)$ and canonical
process $(Y_u, u\ge0)$ with associated canonical filtration
$({\mathcal C}_u^Y, u\ge0)$. Set
\[
{\widetilde{\mathcal G}}_{0}= \bigl(C\bigl(\R^+,\R\bigr),\bigl( {
\mathcal C}_u^Y\bigr), ({\widetilde Q}_{b}, b
\in{\mathcal F}_K) \bigr).
\]
For $A>0$ a $({\mathcal C}_u^Y)$-stopping time, define
the experiment
\[
{\widetilde{\mathcal G}}_{0}^{A}= \bigl(C\bigl(\R^+,\R
\bigr), {\mathcal C}_A^Y, ({\widetilde
Q}_{b}/_{{\mathcal C}_A^Y}, b \in{\mathcal F}_K) \bigr).
\]
For $y \in C(\R^+,\R)$, set ${\widebar A}_0(y)=0$, for $t\in
(t_{i-1}, t_i]$,
%
%e13 #&#
\begin{equation}
\label{abart} {\widebar A}_{t}(y)={\widebar A}_{t_{i-1}}(y) +
\sigma ^2\bigl(y\bigl({\widebar A}_{t_{i-1}}(y)\bigr) \bigr) (t -
t_{i-1}).
\end{equation}
Let ${\widebar T}_u(y)=\inf\{t, {\widebar A}_{t}(y)\ge u\}$.
%
%le3.1 #&#
\begin{lemma}\label{lembar} Set\vspace*{1pt} $(y(u)= x({\bar\tau}_u(x)), u
\ge0)$. Then, ${\widebar A}_{t}(y)={ \bar\rho}_t(x)$ and\break
${\widebar T}_u(y)={\bar\tau}_u(x)$. Consequently, for all
$t\ge0$, ${\widebar A}_t(Y)$ is a $( {\mathcal C}_u^Y)$-stopping time.
\end{lemma}
Thus, the drift term in Proposition~\ref{qtildeb} is ${\bar f}(v)= {\bar f}(v, {\bar\zeta})$ with
\[
{\bar f}(v, y)=\sum_{i\ge1} f\bigl(y\bigl({\widebar
A}_{t_{i-1}}(y)\bigr) \bigr) 1_{({\widebar A}_{t_{i-1}}(y),
{\widebar A}_{t_{i}}(y)]}(v).
\]
%
%%%%%%%%%%%%%%%%%%%%%%%%%%%%%%%%%%%%%%%%%%%%%%%%%%%%%%%
The following result parallel of Proposition~\ref{propdeltaeetilde}
contains the main difficulties.

%th3.1 #&#
\begin{theo} \label{propdeltaggtilde} Assume \textup{(H1)} and \textup{(C)}. For
deterministic $a>0$ and $T=nh$, $
\Delta({\mathcal G}_{0}^{T},{\widetilde{\mathcal G}}_{0}^{{\widebar A}_{T}(Y)})=0
$
and
$
\Delta({\mathcal G}_{0}^{{\bar\tau}_a(X)},{\widetilde{\mathcal
G}}_{0}^{a})=0$.
\end{theo}

The proof uses the following devices.
If $x$ and $y$ are linked by $(x(t), t \ge0) \rightarrow(y(u)=
x({\bar\tau}_u(x)), u\ge0)$, then\vspace*{2pt} $\{u,{\bar\tau
}_u(x)\leq T\} =\{u, u \leq{\widebar A}_T(y)\}$. Similarly, for
$(x(t)= y({\widebar A}_t(y)), t\ge0)$, then $\{ t,{\widebar A}_t(y)
\leq a\}= \{t, t\leq{\bar\tau}_{a}(x)\} $.

%s3.3 #&#
\subsection{Asymptotic equivalence of randomly stopped experiments}\label{sec3.3}

At this point, the triangle inequality implies that, for fixed $T,n,h$
such that $T=nh$,
\[
\Delta\bigl({\mathcal E}_0^T, {\mathcal
G}^{h,n}\bigr)\le\Delta\bigl({\mathcal E}_0^T,
{\mathcal G}_0^{T}\bigr)\le\Delta\bigl({\widetilde{
\mathcal E}}_{0}^{A_T(Y)},{\widetilde{\mathcal
G}}_{0}^{{\widebar A}_T(Y)}\bigr). %
\]
We now introduce the asymptotic framework. Set $T_n=T=nh_n$ and
consider the stopping times
%
%e14 #&#
\begin{equation}
\label{tautaubars} A_n=A_{nh_n}(Y),\qquad{\widebar
A}_n={\widebar A}_{nh_n}(Y), \qquad S_n= {
\widebar A}_n\wedge A_n.
\end{equation}
It remains to study $\Delta({\widetilde{\mathcal E}}_{0}^{A_n},{\widetilde
{\mathcal G}}_{0}^{{\widebar A}_n}) $.
These two experiments have the same sample space but are observed up to
distinct stopping times.
%
%le3.2 #&#
\begin{lemma}\label{rhorhobar} Assume\vspace*{2pt} \textup{(H1)} and \textup{(C)}. There exists a
constant $D$ depending only on $K, K_{\sigma}, \sigma_0, \sigma_1$
such that
$
E_{{\widetilde P}_b} |A_n - {\widebar{A}_n}|\le D n h_n^2$.
\end{lemma}
Using~(\ref{tautaubars}),
the triangle inequality yields
%
%e15 #&#
\begin{equation}
\label{lastineq} \Delta\bigl({\widetilde{\mathcal E}}_{0}^{A_n},{
\widetilde{\mathcal G}}_{0}^{{\widebar A}_n}\bigr) \le \Delta\bigl({
\widetilde{\mathcal E}}_{0}^{A_n},{\widetilde{\mathcal
E}}_{0}^{S_n}\bigr)+ \Delta\bigl({\widetilde{\mathcal
E}}_{0}^{S_n},{\widetilde {\mathcal E}}_{0}^{{\widebar A}_n}
\bigr) + \Delta\bigl({\widetilde{\mathcal E}}_{0}^{{\widebar A}_n},{
\widetilde{\mathcal G}}_{0}^{{\widebar A}_n}\bigr).
\end{equation}
Therefore, we have to study the Le Cam distances, respectively, for the
same experiment observed\vspace*{2pt} up to two distinct times and for two
experiments observed up to
the random time ${\widebar A}_n$. The following holds.
%
%pr3.5 #&#
\begin{prop}\label{proprandom} Assume \textup{(H1)} and \textup{(C)}. There exist
constants $K_1,K_2$ depending only on $K, K_{\sigma}, \sigma_0,
\sigma_1$ such that
%
%e16 #&#
\begin{eqnarray}
\label{propsuffi} \Delta\bigl({\widetilde{\mathcal E}}_{0}^{A_n},{
\widetilde{\mathcal E}}_{0}^{S_n}\bigr)+ \Delta\bigl({
\widetilde{\mathcal E}}_{0}^{{\widebar A}_n},{\widetilde{\mathcal
E}}_{0}^{S_n}\bigr)&\le& K_1 \bigl(
nh_n^2\bigr)^{1/2},
\\
\label{DistanceTV} \Delta\bigl({\widetilde{\mathcal E}}_{0}^{{\widebar A}_n},
{\widetilde{\mathcal G}}_{0}^{{\widebar A}_n}\bigr) &\le& K_2
\bigl(nh_n^2\bigr)^{1/2}.
\end{eqnarray}
Therefore, if $nh_n^2 $ goes to $0$ as $n$ tends to infinity, $\Delta
({\widetilde{\mathcal E}}_{0}^{A_n}, {\widetilde{\mathcal G}}_{0}^{{\widebar
A}_n}) \rightarrow0$.
\end{prop}
Joining Propositions~\ref{propdeltaeetilde},~\ref
{lemeulerdisccont}, Theorem~\ref{propdeltaggtilde} and Proposition
\ref{proprandom} completes the proof of Theorem~\ref{thmainresult}.

%%%%%%%%%%%%%%%%%%%%%%%%%%%%%%%%%%%%%%%%%%%%%%%%%%%%%%%%%%%%%??
%s4 #&#
\section{Concluding remarks}\label{sec4}\label{Concluding}

In this paper, we have obtained the asymptotic equivalence of the
continuous time diffusion~(\ref{eq1}) observed on the time interval
$[0,T]$ and~(\ref{eulerscheme}) the corresponding Euler scheme with
sampling interval $h$ and $T=nh$ in the case of a nonconstant diffusion
coefficient. The discrete Euler scheme model is often used in
applications instead of the diffusion itself. It is broadly accepted as
an appropriate substitute to the diffusion because of its weak
convergence to the diffusion. The equivalence result obtained here was
known for a constant diffusion coefficient. Our contribution is the
extension to the case of a nonconstant diffusion coefficient by means
of random time changed experiments. The constant $K$ in the definition
of the class ${\mathcal F}_K$ is not used for building the Markov
kernels contrary to the diffusion coefficient $\sigma(\cdot)$.
The asymptotic framework is $n\rightarrow+\infty$,
$h=h_n\rightarrow0$ and $nh_n^2=T^2/n\rightarrow0$. In our result,
$T=nh_n$ may be fixed or tend to infinity. We have no assumption
concerning the existence of a stationary regime for the diffusion or
for the Euler scheme. This comes from the assumption that $b$ is
bounded which allows to substantially shorten proofs. For unbounded
drift functions, the two cases ``$T$ bounded'' and ``$T$ tending to
infinity'' have to be distinguished. In the latter case, the diffusion
model must be positive recurrent with moment assumptions on the
stationary distribution.

Compared with other equivalence results, the regularity assumption for
$b$ might seem too strong. However, a classical assumption for
existence and uniqueness of a strong\vspace*{1pt} solution to~(\ref{eq1}) is $b$
locally Lipschitz with linear growth. Generally, authors assume that
$b$ is $C^1$ with linear growth. \citet{DR} consider a special class of
drift functions: $b$ is locally Lipschitz, \textit{known} outside a
compact interval $I$, and H\"older with exponent $\alpha\in(0,1)$
\textit{inside} $I$. They obtain a global asymptotic equivalence of a
stationary diffusion and a mixed Gaussian experiment as $T \rightarrow
+\infty$.

An interesting issue concerns multidimensional diffusions and their
associated Euler scheme. If the diffusion matrix is constant, the
problem is solved [\citet{DR2}]. Otherwise, consider a $d$-dimensional process
$
d\xi_t=b(\xi_t) \,dt + \Sigma(\xi_t) \,dW_t$,
where $b\dvtx  {\mathbb R}^d \rightarrow{\mathbb R}^d$, $\Sigma\dvtx  {\mathbb
R}^d \rightarrow{\mathbb R}^d \otimes{\mathbb R}^d$, $(W_t)$ is a
$d$-dimensional Brownian motion. If $\Sigma(x)$ has the special form
$\Sigma(x)=\sigma(x) P(x)$ where $\sigma\dvtx  {\mathbb R}^d \rightarrow
(0, +\infty)$ and the $d\times d$-matrix $P(x)$ satisfies, for all
$x\in{\mathbb R}^d$, $P(x)P(x)^{t}=I$, the equivalence result is
obtained similarly. Indeed, setting $\rho_t(x)= \int_0^t \sigma
^2(x(s))\,ds$ with inverse $\tau_u(x)$, the time changed\vspace*{2pt} process $\zeta
_u= \xi_{\tau_u(\xi)}$ has a diffusion matrix equal to the identity
matrix and a drift equal to $b(u)/\sigma^2(u)$. As for the continuous
Euler scheme, we can define analogously ${\bar\rho}_t(x)$ and ${\bar\tau}_u(x)$.

Statistical procedures for estimating the drift generally do not use
the knowledge of the diffusion coefficient which appears as a nuisance
parameter. It is an open question to know whether the equivalence
proved here holds when the diffusion coefficient is unknown.

%s5 #&#
\section{Proofs}\label{sec5} \label{Proofs}
\mbox{}
% It relies on classical tools (see {\em e.g.} \citet{KA}, 3.4 and 5.5).
%We have $\zeta_0= \xi_0= \eta$ and
%The change of variable $s=\tau_v(\xi) \Leftrightarrow v=\rho_s(\xi)$
%yields that $dv=\sigma^2 (\xi_s)\,ds=\sigma^2(\zeta_v)\,d\tau_v(\xi)$. So,
%for $0\le v \le a$, $d\tau_v(\xi)=dv/\sigma^2(\zeta_v)$ and $\zeta_u=
%B_u= \int_0^{\tau_u(\xi)}\sigma(\xi_s)\,dW_s.$ We have $<B>_u= \int_0^{
%Brownian motion with respect to the filtration
% $({\mathcal G}_u= {\mathcal A}_{\tau_u(\xi)})$. Note that $({\mathcal
%G}_u)$ satisfies the usual conditions by the continuity of $\tau_.(
%$f=b/\sigma^2$, diffusion coefficient equal to $1$
%and initial condition $ \zeta_0=\eta$, which is ${\mathcal G}_0$
%measurable.
%$\Box$
\begin{pf*}{Proof of Proposition~\ref{propdeltaeetilde}}
By Proposition~\ref{changtemps}, ${\widetilde{\mathcal E}}_{0}^{a}$ is
the image of $ {\mathcal E}_{0}^{\tau_a(X)}$
by the measurable mapping $( x(t), t\in[0,\tau_a(x)]) \rightarrow
(y(u)=x(\tau_u(x)), u\in[0,a])$, which implies $\delta({\mathcal
E}_{0}^{\tau_a(X)},{\widetilde{\mathcal E}}_{0}^{a} ) =0$.

Now, we look at $ {\mathcal E}_{0}^{T} $. As $T=\tau_{\rho_T(x)}(x)$
[(\ref{tau})], the image of $(x(t), t\le T)$ is $ (y(u)= x(\tau
_u(x)), u\le\rho_T(x))$. We must express $\rho_T(x)$ in terms of the
path $y$ and prove that $\rho_T(x)=A_T(y)$. Since $(\rho_T(X)\ge
u)=(\tau_u(X)\le T)$,
%shows that
$\rho_T(X)$ is a stopping time of $({\mathcal C}^X_{\tau_u(X)}, u\ge
0)$. The continuity of $u\rightarrow\tau_u(X)$ implies
%that the following equality of $\sigma$-fields holds:
%
\[
\sigma(X_{\tau_v(X)}, v\le u)= \sigma\bigl(X_s, s\le
\tau_u(X)\bigr). %
\]
Thus, $\rho_T(X)$ is a stopping time of $\sigma(Y_v, v\le u)$ with
$Y_v=X_{\tau_v(X)}$. Observe that, using the change of variable $\tau
_v(X)=s \Leftrightarrow v=\rho_s(X)$, we have
\[
T_u(Y)= \int_0^u \bigl(dv/
\sigma^2(Y_v)\bigr) \,dv = \int_0^{\tau_u(X)}ds= \tau_u(X). %
\]
This implies $\rho_T(X)=A_T(Y)$ which yields $\delta({\mathcal
E}_{0}^{T},{\widetilde{\mathcal E}}_{0}^{A_T(Y)} )=0$.

Consider now the reverse operation. Let $(B_u, u\ge0)$ be a standard
Brownian motion with respect to a filtration $({\mathcal G}_u)$
satisfying the usual conditions and $\zeta_0$ be a ${\mathcal
G}_0$-measurable random variable. We define, for $u\ge0$,
%
%e18 #&#
\begin{equation}
\label{defYub} \zeta_u =\zeta_0 + \int
_0^u \frac{ b(\zeta_v)}{\sigma^2(\zeta
_v)}\,dv +B_u\quad\mbox{and}\quad T_u=T_u(\zeta)=\int_0^u
\frac
{dv}{\sigma^2(\zeta_v)}.
\end{equation}
Clearly, the mapping $u \rightarrow T_u $ is a bijection from $ [0,a]$
onto $[0,T_a]$ with inverse $t\rightarrow T^{-1}(t):=A_t(\zeta)$.
Therefore, we can define, for $0\le t \le T_a$, the process $\xi
_t=\zeta_{A_t(\zeta)}$. The change of variable $
v= A_s(\zeta) \Leftrightarrow s=T_v$ yields that $ ds=dv/ \sigma
^2(\zeta_v)= dv/\sigma^2(\zeta_{A_s(\zeta)})=dv/\sigma^2(\xi_s) $
%$d(\tau^{-1})(s)= \sigma^2(Y_u)= \sigma^2(Y_{\tau^{-1}(s)})$ = X_0+
and equation~(\ref{defYub}) becomes
\[
\xi_t= \xi_0+ \int_0^{A_t(\zeta)}
\frac{b(\zeta_v)}{\sigma^2(\zeta_v)} \,dv+ B_{A_t(\zeta)}= \xi_0+\int
_0^t b(\xi_s) \,ds+
B_{A_t(\zeta)}.
\]
Now, $(M_t=B_{A_t(\zeta)})$ is a martingale w.r.t.
$({\mathcal G}_{A_t(\zeta)})$ satisfying
\[
\langle M\rangle_t= A_t(\zeta) =\int
_0^{A_t(\zeta)}ds= \int_0^t
\sigma^2(\zeta _{A_s(\zeta)})\,ds=\int_0^t
\sigma^2(\xi_s) \,ds.
\]
Hence, $\tau_u(\xi)=A_.(\zeta)^{-1}(u)=T_u$ and $(\xi_t)$ has
distribution $ P_b$. As $(A_t(\zeta))$ is continuous, $({\mathcal
G}_{A_t(\zeta)})$ inherits the usual conditions from $({\mathcal
G}_u)$.

Finally, we can express the above properties on the canonical space.
Let $y=(y(v), v\ge0)$, set $T_u(y)=\int_0^{u} dv/\sigma^2(y(v))$ with
inverse $A_.(y)$ and consider
%the measurable mapping
%
\[
\Psi\dvtx  y \in C\bigl({\mathbb R}^+,{\mathbb R}\bigr)\rightarrow\bigl(x:=y
\bigl(A_t(y)\bigr), t\ge0\bigr)\in C\bigl({\mathbb R}^+,{\mathbb R}
\bigr). %
\]
As $A_t(y)=\int_0^{t}\sigma^2(x(s)) \,ds=\rho_t(x)$, we see that
$A_.(y)^{-1}(u)=\tau_u(x)$. Thus,
$(X_t,t \le\tau_a(X))$ is the image of $(Y(u), u\le a)$ by the
measurable mapping $\Psi$. Hence, $\delta({\widetilde{\mathcal
E}}_{0}^{a},{\mathcal E}_{0}^{\tau_a(X)} ) =0$. Analogously, $(X_t,t
\le T)$ is the image of $(Y(u), u\le A_T(Y))$ which implies $\delta
({\widetilde{\mathcal E}}_{0}^{A_T(Y)},{\mathcal E}_{0}^{T})$.
\end{pf*}

\begin{pf*}{Proof of Proposition~\ref{lemeulerdisccont}}
This proof relies on Lemma~\ref{Lemeuler} below. Define the linear
interpolation between the points $((t_i,Z_i), i\ge0)$:
%
%e19 #&#
\begin{equation}
\label{interp} y(t)= Z_{i}+ \frac{t-t_{i}} {
t_{i+1}-t_{i}} (Z_{i+1}-Z_{i})\qquad\mbox{if } t\in[t_{i}, t_{i+1}]\mbox{ and } i
\ge0.
\end{equation}

%le5.1 #&#
\begin{lemma} \label{Lemeuler}
The solution $(\bar\xi_t)$ of~(\ref{defXbar}) satisfies $
(\bar\xi_{t_i}, i \ge0)=(Z_i, i\ge0)$ where $(Z_i, i\ge0)$ is
the discrete Euler scheme~(\ref{eulerscheme}). Moreover,
%
%e20 #&#
\begin{equation}
\label{BrBr} {\bar\xi}_t= y(t) + \sigma(Z_{i})
B_{i}(t)\qquad\mbox{if } t\in[t_{i}, t_{i+1}]\mbox{ and } i\ge0,
\end{equation}
where\vspace*{2pt} $
B_{i}(t)= W_t-W_{t_{i}} - \frac{t-t_{i}} {t_{i+1}-t_{i}}
(W_{t_{i+1}} - W_{t_i})$. The process $({\bar\xi}_t)$ is adapted
to $({\mathcal A}_t)$,
$((B_i(t),t\in[t_{i}, t_{i+1}]), i\ge0)$ are independent Brownian
bridges and the sequence
$ ((B_i(t),t\in[t_{i}, t_{i+1}]), i\ge0)$ is independent of
$(Z_j, j \ge0)$.
\end{lemma}

This\vspace*{1pt} is a classical result obtained with standard tools.
%{\mathbf Proof of Lemma~\ref{Lemeuler}. }
%
%First, $\bar\xi_0=Z_0=\eta$. By induction, assume that $
%Then
%$$
% - t_i) +\sigma(\bar\xi_{t_i})
%(W_{t_{i+1}} -W_{t_i})=Z_{i+1}
%$$
%Thus, $(\bar\xi_{t_i}, i\ge0) = ( Z_i, i\ge0)$. For $t \in
%[t_{i}, t_{i+1}]$,
%$$
%(W_t-W_{t_{i}} ).
%$$
%Using~(\ref{interp}), $\bar\xi_t =y(t) + \sigma(Z_{i})
%B_{i}(t),$ where $(B_{i}(t)) $ is the Brownian bridge defined in (
%$t_{i} \le t \le t_{i+1}$. The process $(B_{i}(t_i +u), u\in[0,h]) $
%has the distribution of $(W_u -\frac{u}{h}
%W_{h}, 0\le u \le h)$.
%Using that, for all $i \ge0$, $B_i(t)$ is ${\mathcal A}_{t_{i+1}}$-
%measurable and independent of ${\mathcal A}_{t_{i}}$ yields that
%$(B_i, i\ge0)$
%are independent processes. Elementary computations yield that
%$(B_i(\cdot), i\ge0)$ is independent of the vector $(W_{t_{i+1}}
%-W_{t_i}, i\ge0)$. Since $Z_0$ is independent of $ (W_t, t\ge0)$, we
%first get the independence of $Z_0$ and $(B_i(\cdot),W_{t_{i+1}} -W_{t_i},
%i\ge0)$. This implies that $Z_0$, $(B_i(\cdot), i\ge0)$, $(W_{t_{i+1}}
%-W_{t_i}, i\ge0)$ are independent. As $\sigma(Z_i, i\ge0) \subset
%
We may now complete the proof of Proposition~\ref{lemeulerdisccont}.
Since
%the Euler scheme
$(Z_i, i\ge0)$ is the image of\vspace*{1pt}
%the sample path
$(\bar{\xi}_{t}, t \ge0)$ by the mapping $x(\cdot) \rightarrow
(x(t_i), i\ge0)$, $ \delta({\mathcal G}_{0}^{Nh},{\mathcal
G}^{h, N})=0$.

Consider, for $\omega\in\Omega$,
the application
$
\Phi_{\sigma,h} = \Phi\dvtx \R^{\N} \rightarrow{\mathcal C} (\R^+, \R)
$
defined by $(x_i, i\ge0) \rightarrow x(\cdot)$ with
$x(t) = x_{i-1}+ \frac{t-t_{i-1}} {
t_i-t_{i-1}} (x_i - x_{i-1}) + \sigma(x_{i-1})B_{i-1} (t,\omega)
$ for $ t \in[t_{i-1}, t_i]$.
%,  i\ge0$.
As $\sigma$ is known, $\Phi$ is a randomization and, by Lemma~\ref
{Lemeuler},
%the experiment
$ {\mathcal G}_{0}^{N h}$ is the image by
$\Phi$ of ${\mathcal G}^{h, N}$. Hence, $ \delta({\mathcal G}^{h, N},
{\mathcal G}_{0}^{Nh})= 0$.
%Hence, the result.
\end{pf*}

\begin{pf*}{Proof of Proposition~\ref{qtildeb}}
By definition of $({\bar\zeta}_u)$, we have
%
%e21 #&#
\begin{equation}
\label{Ybarbis} {\bar\zeta}_u = {\bar\xi}_0+\int
_0^{{\bar\tau
}_u({\bar\xi})} \sum_{i\ge0}
b({\bar\xi}_{t_i})1_{t_i < s\le t_{i+1}} \,ds + {\widebar B}_u,
\end{equation}
where ${\widebar B}_u = \int_0^{{\bar\tau}_u({\bar\xi})} \sum_{i\ge0}
\sigma({\bar\xi}_{t_i}) 1_{t_i < s \le t_{i+1}} \,dW_s$
%So $ ({\widebar B}_u)$
is a martingale w.r.t.
${\widebar{\mathcal G}}_u= {\mathcal A}_{{\bar\tau
}_u({\bar\xi})}$
with quadratic\vspace*{1pt} variations
$ \langle  {\widebar B} \rangle_u= \int_0^{{\bar\tau}_u({\bar\xi})}
\sum_{i\ge0} \sigma^2({\bar\xi}_{t_i}) 1_{t_i < s\le t_{i+1}} \,ds=u$.
Therefore, $({\widebar B}_u)$ is a Brownian motion with respect to
$({\widebar{\mathcal G}}_u)$.

In the integral of~(\ref{Ybarbis}), the change of variable
$ s = {\bar\tau}_v({\bar\xi})\Leftrightarrow v
= {\bar\rho}_s({\bar\xi}) $ yields, noting that
$dv=\sigma^2( {\bar\xi}_{t_i}) \,ds $
for $ v \in({\bar\rho}_{t_i}({\bar\xi}), {\bar\rho
}_{t_{i+1}}
({\bar\xi})]$,
and that $t_i={\bar\tau}_{{\bar\rho}_{t_i}(x)}(x)$,
\begingroup
\abovedisplayskip=9pt
\belowdisplayskip=9pt
%e22 #&#
\begin{equation}
{\bar\zeta}_u = {\bar\xi}_0+\int_0^{u}
\sum_{i\ge
0}\frac{ b({\bar\xi}_{t_i})} {
\sigma^2({\bar\xi}_{t_i})} 1_{{\bar\rho
}_{t_i}({\bar\xi}) < v \le
{\bar\rho}_{t_{i+1}}({\bar\xi})} \,dv +
{\widebar B}_u,
\end{equation}
where ${\bar\xi}_{t_i}= {\bar\zeta}_{{\bar\rho
}_{t_i}({\bar\xi})}=Z_i$
is the discrete Euler scheme (Lemma~\ref{Lemeuler}).

Thus,\vspace*{1pt} $({\widebar Y}_u)$ defined in~(\ref{Ybar}) is a process with
diffusion coefficient equal to $1$ and drift term ${\bar f}(v)$.
We now\vspace*{2pt} check that ${\bar f}(v)$ is predictable w.r.t.
$({\widebar{\mathcal G}}_u)$, that is, $ \forall i$,
${\bar\rho}_{t_i}
({\bar\xi})$ is a $({\widebar{\mathcal G}}_u)$-stopping time\vspace*{1pt} and
${\bar\zeta}_{{\bar\rho}_{t_i}({\bar\xi})}$ is
${\widebar{\mathcal G}}_{{\bar\rho}_{t_i}({\bar\xi})}
$-measurable.
Noting that
% for all $u$,
$ ({\bar\rho}_{t_i}({\bar\xi})\le u)=
({\bar\tau}_u({\bar\xi})
\ge t_i)
$
belongs to ${\widebar{\mathcal G}}_u= {\mathcal A}_{{\bar\tau
}_u({\bar\xi})}$
yields that ${\bar\rho}_{t_i}({\bar\xi})$ is a
\mbox{$({\widebar{\mathcal G}}_u)$-}stopping time. We know that
$
{\bar\zeta}_{{\bar\rho}_{t_i}({\bar\xi
})}={\bar\xi}_{t_i}
$
is $ {\mathcal A}_{t_i}$-measurable, which achieves the proof since
\mbox{${\mathcal A}_{t_i}=
{\widebar{\mathcal G}}_{{\bar\rho}_{t_i}({\bar\xi})}$.}
\end{pf*}

\begin{pf*}{Proof of Lemma~\ref{lembar}}
The relation $(y(u)=x({\bar\tau}_u(x))$ is equivalent to
$(y({\bar\rho}_t(x))= x(t))$.
First, note that ${\widebar A}_{t_1}(y)= \sigma^2(y(0))t_1= \sigma
^2(x(0))t_1= {\bar\rho}_{t_1}(x)$. By induction, assume that
${\widebar A}_{t_j}(y)={\bar\rho}_{t_j}(x)$ for $j\le i-1$. Then
\begin{eqnarray*}
{\widebar A}_{t_i}(y) &=& {\bar\rho}_{t_{i-1}}(x)+ \sigma
^2\bigl(y\bigl({\bar\rho}_{t_{i-1}}(x)\bigr)\bigr)
(t_i- t_{i-1})
\\
&=&  {\bar\rho }_{t_{i-1}}(x)+
\sigma^2\bigl(x(t_{i-1})\bigr) (t_i-
t_{i-1})
= {\bar\rho }_{t_{i}}(x).
\end{eqnarray*}
Thus, the two inverse functions coincide: ${\widebar T}_u(y)=
{\bar\tau}_u(x)$. As above, we deduce that $A_t(y)$ is a
stopping time w.r.t. $({\mathcal C}_u^Y)$ with $Y_u=X_{{\bar\tau}_u(X)}$.
\end{pf*}\endgroup

\begin{pf*}{Proof of Theorem~\ref{propdeltaggtilde}}
The\vspace*{1pt} proof is divided in several steps.

First, as
${\widetilde{\mathcal G}}_0^{a}$ is the image of
${\mathcal G}_0^{{\bar\tau}_a(X)}$ by the measurable mapping
$(x(t), t\le{\bar\tau}_a(x))\rightarrow(y(u)=x({\bar\tau
}_u(x)), u \in[0,a])$,
$\delta({\mathcal G}_0^{{\bar\tau}_a(X)},{\widetilde{\mathcal
G}}_0^{a})=0$.

Now consider ${\mathcal G}_0^{T}$. We have $T= {\bar\tau
}_{{\bar\rho}_T(x)}(x)$. Hence, the image of $(x(t), t\le T)$ is
$(y(u)=x({\bar\tau}_u(x)), u\le{\bar\rho
}_T(x)={\widebar A}_T(y))$ according to Lemma~\ref{lembar}. This
proves that $\delta({\mathcal G}_0^{T},{\widetilde{\mathcal
G}}_0^{{\widebar A}_T(y)})=0$.

Let us study the other deficiencies.
%prove that $\delta({\widetilde{\mathcal G}}_0^{a},{\mathcal G}_0^{{
We first construct a process $({\bar\zeta}_u)$ with distribution
${\widetilde Q}_b$ (step~1), then
a process $({\bar\xi}_t)$ with\vspace*{1pt} distribution $ Q_b$ obtained from
$({\bar\zeta}_u)$ by the mapping
$(y(u), u \ge0) \rightarrow(y({\widebar A}_t (y)),  t \ge0)$ (step~2).

%{\bar\tau}_a(x))\rightarrow(x({\bar\tau}_u(x)), u \in
%[0,a])$,
%$\delta({\mathcal G}_0^{{\bar\tau}_a(x)},{\widetilde{\mathcal
%G}}_0^{a})=0$.\\

\begin{longlist}[\textit{Step} 1.]
\item[\textit{Step} 1.]
%Let us consider now the reverse transformation.
Let $( {\widebar B}_u)$ be a Brownian motion w.r.t. a filtration
$( {\widebar{\mathcal G}}_u)$ satisfying the usual conditions. Assume
that ${\bar\zeta}_0 $ is
$ {\widebar{\mathcal G}}_0$-measurable.
Then we define recursively a sequence of random times $(T_i)$ and a continuous
process $( {\bar\zeta}_u)$. First, set
$ T_0=0$, then
\begin{eqnarray}
\qquad T_1 &=& T_1({\bar\zeta})=\sigma^2({\bar
\zeta}_0)t_1, \qquad {\bar\zeta}_u= {\bar
\zeta}_0+ f({\bar\zeta}_0) u+ {\widebar
B}_u\qquad\mbox{for } 0<u\le T_1, \nonumber
\\
%
%e23 #&#
\label{defTi} T_i &=& T_i({\bar\zeta}) =T_{i-1}
+ \sigma^2({\bar\zeta}_{T_{i-1}}) (t_i-t_{i-1}),
\\
%
%e24 #&#
\label{defYTi} {\bar\zeta}_u &=& {\bar\zeta}_{T_{i-1}} +f( {\bar
\zeta}_{T_{i-1}}) (u-T_{i-1})+{\widebar B}_u-{\widebar
B}_{T_{i-1}}\qquad\mbox{for }T_{i-1} <u \le T_i.
\end{eqnarray}
Note that $T_i={\widebar A}_{t_i}(\zeta)$ [see~(\ref{abart})].
%%%%%%%%%%%%%%%%%%%%%%%%%%%%%%%%%%%%%%%%%%%%%%%%%%%%%%%%%%%%%%%%%%%%%%
%
%le5.2 #&#
\begin{lemma} \label{lemTiYTi}
The\vspace*{1pt} sequence ($T_i$) is an increasing sequence of
$ ({\widebar{\mathcal G}}_u)$-stopping times such that, for all $i\ge
1$, $ T_i$ is
$ {\widebar{\mathcal G}}_{T_{i-1}}$ measurable. Moreover, the process
$ ({\bar\zeta}_{u}) $ defined in~(\ref{defTi}), (\ref
{defYTi}) is a diffusion-type process adapted to $ ({\widebar{\mathcal G}}_u)$
with diffusion coefficient equal to $1$ and drift coefficient
\[
{\bar f}(u, y)= \sum_{i\ge1} f\bigl(y
\bigl(T_{i-1}(y)\bigr)\bigr) 1_{T_{i-1}(y)< u
\le T_i(y)}, %
\]
where\vspace*{1pt} $(T_i(y)={\widebar A}_{t_i}(y), i\ge0)$ are recursively defined
as in~(\ref{abart}) using $y(\cdot)$ and $f=b/\sigma^2$ [see~(\ref{f})].
Hence, the process $({\bar\zeta}_{u})$ has
distribution ${\widetilde Q}_b$.
\end{lemma}
%
%%%%%%%%%%%%%%%%%%%%%%%%%%%%%%%%%%%%%%%%%%%%%%%%%%%%%%%%%%%%%%%%%%%%%%%%%%%
%
\begin{pf}
First,\vspace*{1pt} $T_1$ is ${\widebar{\mathcal G}}_0$-measurable, thus $\{T_1
\leq u\}
\in{\widebar{\mathcal G}}_0
\subset{\widebar{\mathcal G}}_u$. Hence, $T_1$ is a $({\widebar
{\mathcal G}}_u)$-stopping time. Now,\vspace*{1pt} ${\bar\zeta}_{u}=
{\bar\zeta}_{0} +f ( {\bar\zeta}_{0}) u+ {\widebar
B}_{u} $ is
%for all $u$,
$ {\widebar{\mathcal G}}_{u} $-measurable. Thus, $T_1$ and
${\bar\zeta}_{T_1}$ are
$ {\widebar{\mathcal G}}_{T_{1}}$ measurable.

By induction, assume that, for $ 1\le j \le i$, $T_j$ is ${\widebar
{\mathcal
G}}_{T_{j-1}}$-measurable,
$T_j$ is a \mbox{$({\widebar{\mathcal G}}_u)$-}stopping time,\vspace*{2pt} and
$({\bar\zeta}_u,
%$ defined by~(\ref{defYTi}) for $
u\le T_i)$ is ${\widebar{\mathcal G}}_{u}$-measurable.
Now, for $u>T_i$,
${\bar\zeta}_u ={\bar\zeta}_{T_{i}}
+f ( {\bar\zeta}_{T_{i}})(u-T_{i})+{\widebar B}_u-{\widebar
B}_{T_{i}}$ defined by~(\ref{defYTi}) is ${\widebar{\mathcal
G}}_{u}$-measurable.
As $T_{i+1}=T_{i} + \sigma^2({\bar\zeta
}_{T_{i}})(t_{i+1}-t_{i})$, the induction assumption
yields that $T_{i+1}$ is ${\widebar{\mathcal G}}_{T_{i}}$-measurable
and, since $T_i<T_{i+1} $ by~\textup{(C)},
%and using the induction assumption,
%
\[
\forall v \ge u\qquad\{T_{i+1} \le u\}= \{T_{i+1} \le u\}
\cap\{T_{i} \le v \} \in{\widebar {\mathcal G}}_v.
\]
This\vspace*{1pt} implies that $\{T_{i+1} \le u\} = \{T_{i+1} \le u\}\cap\bigcap_{v > u} \{T_{i} \le v\} \in
\bigcap_{v > u} {\widebar{\mathcal G}}_v= {\widebar{\mathcal G}}_u$
which proves that $T_{i+1}$ is a
$({\widebar{\mathcal G}}_u)$-stopping time. Thus, $T_{i+1}$ and
${\bar\zeta}_{T_{i+1}}$ are
\mbox{${\widebar{\mathcal G}}_{T_{i+1}}$-}measurable.
The proof of Lemma~\ref{lemTiYTi} is now complete.
\end{pf}

\item[\textit{Step} 2.]
Let us study the distribution of ${\bar\xi}_t$ defined as
%
%e25 #&#
\begin{equation}
\label{xbarbis} {\bar\xi}_t= {\bar\zeta}_{{\widebar A}_t({\bar\zeta})}.
\end{equation}
%
%and look at the distribution of $({\bar\xi}_t) $.
By Lemma~\ref{lemTiYTi},
${\widebar A}_{t_i}({\bar\zeta})= T_i$ is a $({\widebar
{\mathcal G}}_u)$-stopping time.
%We check that, $\forall t,   {\widebar A}_{t}({\bar\zeta})$ is
%a
%$({\widebar{\mathcal G}}_u)$ -stopping time.
For $t_i \leq t \leq t_{i+1}$, ${\widebar A}_{t}({\bar\zeta})=
T_i + (t - t_i) \sigma^{2}({\bar\zeta}_{T_i})$ is
${\widebar{\mathcal G}}_{T_i}$-measurable, so
\[
\forall v > u\qquad\bigl\{ {\widebar A}_{t}({\bar\zeta}) \le u\bigr\}
=\bigl\{ {\widebar A}_{t}({\bar\zeta}) \le u\bigr\} \cap
\{T_{i} \le v \} \in{\widebar {\mathcal G}}_v.
\]
Hence,\vspace*{1pt} $\{ {\widebar A}_{t}({\bar\zeta}) \le u\} =\{ {\widebar
A}_{t}({\bar\zeta}) \le u\} \cap\bigcap_{v > u} \{T_{i} \le
v\} \in
\bigcap_{v > u} {\widebar{\mathcal G}}_v= {\widebar{\mathcal G}}_u$
which proves that ${\widebar A}_{t}({\bar\zeta}) $ is a
$({\widebar{\mathcal G}}_u)$-stopping time.

Thus, we can define the filtration $({\widebar{\mathcal A}}_t:=
{\widebar{\mathcal G}}_{ {\widebar A}_{t}({\bar\zeta}) })$ to
which $({\bar\xi}_t)$ is adapted.
%%%%%%%%%%%%%%%%%%%%%%%%%%%%%%%%%%%%%%%%%%%%%%%%%%%%%%%%%%%%%%%%%%%%%%%%%%%%%%%%%%%%%%%%
%%%%%%%%%%%%%%%%%%%%%%%%%%%%%%%%%%%%%%%%%%%%%%%%%%%%%%%%%%%%%%%%%%%%%%%%%%%%%%%%%%%%%%
%By construction, ${\bar\tau}_a({\widebar Y}_.) $ is a $({

%le5.3 #&#
\begin{lemma} \label{eulerretour}
The sequence $({\bar\xi}_{t_i}={\bar\zeta}_{T_{i}}, i\ge
0)$, with $({\bar\zeta}_u)$ defined by
(\ref{defTi})--(\ref{defYTi}), $({\bar\xi}_t)$ in (\ref
{xbarbis}), has the distribution of the discrete Euler scheme (\ref
{eulerscheme}).
\end{lemma}
\begin{pf}
For all $i\ge0$, the process
%
%e26 #&#
\begin{equation}
\label{bbari} \bigl({\widebar B}^{(i)}_v={\widebar
B}_{T_{i}+v}-{\widebar B}_{T_{i}}, v\ge0\bigr)
\end{equation}
is a Brownian motion independent of ${\widebar{\mathcal
G}}_{T_i}={\mathcal A}_{t_i}$, adapted to
$({\widebar{\mathcal G}}_{T_i+v})$. As ${\bar\xi
}_{t_i}={\bar\zeta}_{T_{i}}$ is ${\widebar{\mathcal
G}}_{T_i}$-measurable,
this r.v. is independent of $({\widebar B}^{(i)}_v, v\ge0)$. Define
%
%e27 #&#
\begin{equation}
\label{eps} \varepsilon_{i+1}= \frac{{\widebar B}_{T_{i+1}}-{\widebar
B}_{T_i}}{\sqrt{T_{i+1}-T_i}}=
\frac{{\widebar B}^{(i)}_{ \sigma^{2}({\widebar
Y}_{T_i})(t_{i+1}-t_i)}}{ \sigma({\widebar Y}_{T_i})\sqrt{t_{i+1}-t_i}}. %\frac{{\widebar B}_{T_{i}+\sigma^{2}({\widebar
%Y}_{T_i})(t_{i+1}-t_i)}-
%{\widebar B}_{T_{i}}}{\sigma({\widebar Y}_{T_i})\sqrt{t_{i+1}-t_i}}
\end{equation}
The random variable $\varepsilon_{i+1}$ is ${\widebar{\mathcal
G}}_{T_{i+1}}$-measurable.
We can write
%
%e28 #&#
\begin{equation}
\label{schemabari} {\bar\zeta}_{T_{i+1}}={\bar\zeta}_{T_{i}}+ b({\bar
\zeta}_{T_{i}}) (t_{i+1}-t_i) + \sigma({\bar
\zeta}_{T_{i}}) \sqrt{t_{i+1}-t_i} \varepsilon
_{i+1}, \qquad i\ge0.
\end{equation}
To\vspace*{2pt} conclude, it is enough to prove that $(\varepsilon_i, i\ge1)$ is a
sequence of
i.i.d. standard Gaussian random variables, independent of
${\widebar{\mathcal G}}_{0}$.

Applying Proposition~\ref{propauxi} of the \hyperref[secappen]{Appendix} yields that, for
all $i\ge0$, $\varepsilon_{i+1}$ is a standard Gaussian variable
independent of ${\widebar{\mathcal G}}_{T_i}$. This holds for $i=0$
and proves that $\varepsilon_1$ is independent of ${\widebar
{\mathcal G}}_{0}$ and has distribution ${\mathcal N}(0, 1)$. By
induction, assume that $(\varepsilon_k, k\le i-1)$ are
i.i.d. standard Gaussian random variables, independent of~${\widebar{\mathcal G}}_{0}$.
Consider ${\bar\zeta}_0\sim\eta$. As $({\bar\zeta}_0,
\varepsilon_k, k\le i-1)$ is ${\widebar{\mathcal
G}}_{T_i}$-measurable, we get that $\varepsilon_{i+1}$ is a standard
Gaussian variable independent of $({\bar\zeta}_0, \varepsilon
_k, k\le i-1)$.
This completes the proof of Lemma~\ref{eulerretour}.
\end{pf}
Define now $({\bar x}(t))$ as the linear interpolation between the
points $(t_i, {\bar\xi}_{t_i})$.
We now describe the processes $({\bar\xi}_t- {\bar x}(t))$
for $t_i \leq t \leq t_{i+1}$.
%%%%%%%%%%%%%%%%%%%%%%%%%%%%%%%%%%%%%%%%%%%%%%%%
%
%le5.4 #&#
\begin{lemma} \label{eulercontinuretour}
For\vspace*{2pt} $t \in[t_i, t_{i+1}]$,
$
{\bar\xi}_{t} = {\bar x}(t)+ \sigma({\bar\xi}_{t_i}
) {\widebar C}_{i}(t)$,
where $(({\widebar C}_{i}(t), t_i \le t \le t_{i+1}), i\ge0)$ is a
sequence of independent Brownian bridges adapted to
$({\widebar{\mathcal A}}_t )$,
%. The sequence $(({\widebar C}_{i}(t), t_i \le t \le t_{i+1}), i\ge0)$ is
independent of $({\bar\xi}_{t_j}, j \ge0)$.
\end{lemma}
\begin{pf} We have
$
{\bar y}(u)= {\bar\zeta}_{T_{i}}+ \frac{u -
T_i}{T_{i+1}-T_i} ({\bar\zeta}_{T_{i+1}}-{\bar\zeta}_{T_{i}})$.
Using~(\ref{eps})--(\ref{schemabari}), we obtain, for $u \in[T_i, T_{i+1}]$,
\begin{eqnarray*}
{\bar\zeta}_{u}& =& {\bar y}(u) + {\widebar B}_{u}-{
\widebar B}_{T_i}- \frac{u - T_i}{T_{i+1}-T_i} \sigma ({\bar\zeta}_{T_{i}})
\sqrt{t_{i+1}-t_i} \frac{{\widebar B}_{T_{i+1}}-{\widebar B}_{T_i}}{\sqrt {T_{i+1}-T_i}}
\\
&=& {\bar y}(u) + D_{i}(u)
\end{eqnarray*}
with
\[
D_{i}(u)= {\widebar B}_{u}-{\widebar B}_{T_i}- \frac{u - T_i}{T_{i+1}-T_i} ({
\widebar B}_{T_{i+1}}- {\widebar B}_{T_i}).
\]
%
%D_{i}(u)= {\widebar B}_{u}-{\widebar B}_{T_i}- \frac{u -
%T_i}{T_{i+1}-T_i} ({\widebar
%B}_{T_{i+1}}- {\widebar B}_{T_i}).
For\vspace*{1pt} $t_i \leq t \leq t_{i+1}$, using~(\ref{abart}) and (\ref
{defTi}), we get
$ {\bar x}(t)= {\bar y} ({\widebar A} _t({\bar\zeta}_.))$.
Thus,
${\bar\xi}_{t}-{\bar x}(t)= D_i({\widebar A} _t({\bar\zeta}_.))$, and
define, using~(\ref{bbari}), ${\widebar C}_{i}(t)$ by
\[
{\bar\xi}_{t}-{\bar x}(t)= {\widebar B}^{(i)}_{\sigma
^2({\bar\xi}_{t_i})(t-t_i)}
- \frac{t - t_i}{t_{i+1}-t_i} {\widebar B}^{(i)}_{\sigma^2({\bar\xi}_{t_i})(t_{i+1}-t_i) }= \sigma({\bar
\xi}_{t_i}) {\widebar C}_{i}(t).
\]
Proving\vspace*{2pt} that $({\bar\xi}_{t_i}, i\ge0)$ is independent of
$(({\widebar C}_{i}(t), t\in[t_i, t_{i+1}]), i \ge0)$ is equivalent to\vspace*{2pt}
proving that $({\bar\xi}_{0},\varepsilon_i, i\ge1)$ is
independent of $(({\widebar C}_{i}(t), t\in[t_i, t_{i+1}]),\break  i \ge0)$. We
now show that, $\forall i\ge1, ({\bar\xi}_{0},\varepsilon_1,
\ldots,\varepsilon_i)$ is independent of $({\widebar
C}_{0},\ldots,\break
{\widebar C}_{i-1})$ and that the latter processes are independent Brownian
bridges. Using Proposition~\ref{propauxi} with $B={\widebar
B}^{(i-1)}$, ${\mathcal F}_.={\widebar{\mathcal G}}_{T_{i-1}+.}$,
$\tau= \sigma^{2}({\bar\xi}_{t_{i-1}})$, $i\ge1$ yields
that $
W_i(t-t_{i-1})= \frac{1}{\sigma({\bar\xi
}_{t_{i-1}})}{\widebar B}^{(i-1)}_{\sigma^{2}({\bar\xi
}_{t_{i-1}})(t-t_{i-1})}$, $
t\ge t_{i-1}$,
%{\sigma({\bar\xi}_{t_{i-1}})},
is a Brownian motion independent of ${\widebar{\mathcal
G}}_{T_{i-1}}$. Thus, $({\widebar C_{i-1}}(t), t\in[t_{i-1}, t_i])$ is a
Brownian\vspace*{2pt} bridge independent of $W_i(t_i-t_{i-1})= \varepsilon_i \sqrt {t_i-t_{i-1}}$. Moreover, ${\widebar{\mathcal G}}_{T_{i-1}}$,
$W_i(t_i-t_{i-1})$, and $({\widebar C_{i-1}}(t), t\in[t_{i-1}, t_i])$\vspace*{2pt} are
independent.

For $i=1$, as ${\bar\xi}_0$ is ${\widebar{\mathcal
G}}_{0}$-measurable, we get that ${\bar\xi}_0$, $\varepsilon
_1$, ${\widebar C}_0$ are independent and ${\widebar C}_0$ is a Brownian
bridge. By induction, let us assume that ${\bar\xi}_0,
\varepsilon_1, \ldots, \varepsilon_i,\break {\widebar C}_0, \ldots, {\widebar
C}_{i-1}$ are independent and that ${\widebar C}_0, \ldots, {\widebar
C}_{i-1}$ are Brownian bridges (on their respective interval of
definition). As $Z=({\bar\xi}_0, \varepsilon_1, \ldots,
\varepsilon_i,{\widebar C}_0, \ldots, {\widebar C}_{i-1})$ is ${\widebar
{\mathcal G}}_{T_{i-1}}$-measurable, we get that $Z, \varepsilon
_{i+1}, {\widebar C}_i$ are independent. The proof of Lemma~\ref
{eulercontinuretour} is complete.
\end{pf}
Thus, we have constructed a process $({\bar\xi}_t)$ with
distribution $Q_b$ obtained by the mapping $(y(u), u \ge0) \rightarrow
(y({\widebar A}_t (y)),  t \ge0)$. Hence, $(x(t)=y({\widebar A}_t
(y)), t \le{\widebar T}_a(y)={\bar\tau}_a(x))$ is the image of
$(y(u), u\le a)$. This proves $\delta({\widetilde{\mathcal
G}}_0^{a},{\mathcal G}_0^{{\bar\tau}_a(X)})=0$.\vspace*{2pt}

Moreover, $(x(t), t\le T)$ is the image of $(y(u), u\le{\widebar A}_T
(y))$. This yields $\delta({\widetilde{\mathcal G}}_0^{{\widebar
A}_T(Y)},{\mathcal G}_0^{T})=0$. This completes the proof of Theorem
\ref{propdeltaggtilde}.\quad\qed
\end{longlist}\noqed
\end{pf*}

\begin{pf*}{Proof of Lemma~\ref{rhorhobar}}
Using~(\ref{tuat}), $ A_t(y)=u \Leftrightarrow T_u(y)=t$ yields that $
A_n=\int_0^{nh_n} \sigma^2(y(A_s(y))) \,ds$. Combining with (\ref
{abart}), we get
\[
A_n - {\widebar A}_n= \sum
_{i=1}^n \int_{t_{i-1}}^{t_{i}}
\bigl( \sigma ^2\bigl(y\bigl(A_s(y)\bigr)\bigr) -
\sigma^2\bigl(y\bigl({\widebar A}_{t_{i-1}}(y)\bigr) \bigr)
\bigr)\,ds. %
\]
Under\vspace*{2pt} ${\widetilde P}_b$, $(Y(A_t(Y))=X_t)$ has distribution $P_b$ (see
proof of Proposition~\ref{propdeltaeetilde}). Hence,
$
E_{ {\widetilde P}_b} |A_n - {\widebar{A_n}}|=E_{P_b}| \sum_{i=1}^n \int_{t_{i-1}}^{t_{i}}  (\sigma^2(X_s)- \sigma^2(X_{t_{i-1}}) )\,ds|
$.
Denoting by ${\mathcal L}$ the generator of the diffusion $(X_t)$
(${\mathcal L}h =(1/2)\sigma^2h'' + b h'$), the It\^o formula yields
$
\int_{t_{i-1}}^{t_{i}}  (\sigma^2(X_s) - \sigma
^2(X_{t_{i-1}}) )\,ds=B_1(i)+B_2(i)$,
with
\begin{eqnarray*}
B_1(i) &=&\int_{t_{i-1}}^{t_{i}}dv\int
_{t_{i-1}}^{s} {\mathcal L}\sigma ^2
(X_u) \,du,
\\
B_2(i) &=& \int_{t_{i-1}}^{t_{i}}dv
\int_{t_{i-1}}^{s} \bigl(\sigma^2
\bigr)'(X_u) \sigma(X_u)\,dB_u.
\end{eqnarray*}
Condition \textup{(C)} and \textup{(H1)} ensure that $\|{\mathcal L}\sigma^2\|_{\infty}
$ is bounded by $D_1$ depending on $K, K_{\sigma}, \sigma_1$, so
that, $
|B_1(i)| \le D_1 h_n^2/2$.
For the second term,
\begin{eqnarray*}
B_2(i)&=& \int_{t_{i-1}}^{t_{i}}ds \int
_{t_{i-1}}^{s} \bigl(\sigma ^2
\bigr)'(X_u)\sigma(X_u) \,dB_u
\\
&=& \int_{t_{i-1}}^{t_{i}}(t_i -u) \bigl(\sigma
^2\bigr)'(X_u) \sigma(X_u)\,dB_u,
\\
\sum_{i=1}^n B_2(i)&=&\int
_0^{nh_n} H_u^{(n)}
\,dB_u,
\end{eqnarray*}
where
\[
H_u^{(n)}= \sum
_{i=1}^{n} 1_{]t_{i-1}, t_i]}(u) (t_i
-u) \bigl(\sigma ^2\bigr)'(X_u)
\sigma(X_u). %
\]
This yields
$
E_{ P_b} (\sum_{i=1}^n B_2(i))^2 = E_{ P_b}\int_0^{nh_n}
(H_u^{(n)})^2 \,du \le D_2   n h_n^3
$ with $D_2$ a constant.
Therefore,
$
E_{ {\widetilde P}_b} |A_n - {\widebar{A}_n}|\le D' (nh_n^2+
(nh_n^3)^{1/2})\le D nh_n^2$.
\end{pf*}

\begin{pf*}{Proof of Proposition~\ref{proprandom}}
\textit{Proof of inequality}~(\ref{propsuffi}).
As ${\widetilde{\mathcal E}}_0^{S_n}$ is a restriction of ${\widetilde
{\mathcal E}}_0^{A_n}$ to a smaller $\sigma$-algebra, $\delta({\widetilde
{\mathcal E}}_0^{A_n},{\widetilde{\mathcal E}}_0^{S_n})=0$. To\vspace*{1pt} evaluate
the other deficiency, we introduce a kernel from ${\widetilde{\mathcal
E}}_0^{S_n}$ to ${\widetilde{\mathcal E}}_0^{A_n}$. Let $B \in{\mathcal
C}_{A_n}^Y$, and set
$
N(\omega, B)= E_{{\widetilde P}_0}(1_B | {\mathcal C}_{S_n}^Y)(\omega)$,
where ${\widetilde P}_0$, corresponding to $b=0$, is the distribution of
$(\eta+B_u, u\ge0)$.
Now, $N({\widetilde P}_b|{\mathcal C}_{S_n}^Y)$ defines\vspace*{1pt} a probability on
$(C(\R^+,\R), {\mathcal C}_{A_n}^Y)$ with density w.r.t. ${\widetilde
P}_0|{\mathcal C}_{A_n}^Y$, $ (d{\widetilde P}_b/d{\widetilde P}_0)|{\mathcal
C}_{S_n}^Y$. Indeed,
for $B \in{\mathcal C}_{A_n}$,
\begin{eqnarray*}
N({\widetilde P}_b|{\mathcal C}_{S_n}) (B)& =& \int
_{\Omega}N(\omega, B)\,d\bigl({\widetilde P}_b|{
\mathcal C}_{S_n}^Y\bigr)= E_{{\widetilde P}_0} \biggl(
\frac{d{\widetilde P}_b}{d{\widetilde P}_0}\bigg|{\mathcal C}_{S_n}^Y E_{{\widetilde
P}_0}
\bigl(1_B | {\mathcal C}_{S_n}^Y\bigr) \biggr)
\\
& =& E_{{\widetilde P}_0} \biggl(\frac{d{\widetilde P}_b}{d{\widetilde
P}_0}\bigg|{\mathcal C}_{S_n}^Y1_B
\biggr).
\end{eqnarray*}
%
%Thus, $N({\widetilde P}_b|{\mathcal C}_{S_n})$ has density w.r.t. ${\widetilde
%P}_0|{\mathcal C}_{A_n}^Y$ equal to
%$ (d{\widetilde P}_b/d{\widetilde P}_0)|{\mathcal C}_{S_n}^Y.$
For $T$ a bounded stopping time,
\[
\frac{d{\widetilde P}_b}{d{\widetilde P}_0}\bigg|{\mathcal C}_{T}^Y={\widetilde
L}_T(b)=\exp{\biggl(\int_0^T
f(Y_u) \,dY_u - \int_0^T
\frac{1}{2} f^2(Y_u) \,du\biggr)}. %
\]
%
%Denote by ${\widetilde L}_T(b)$ the density $\frac{d{\widetilde P}_b}{d{\widetilde
%P}_0}|{\mathcal C}_{T}$. We have
%$$
%{\widetilde L}_T(b)=\exp{(\int_0^T f(X_u) \,dX_u - \int_0^T \frac{1}{2}
%f^2(X_u) \,du)}
%$$
%with $f=b/\sigma^2$.
Thus, $ (d{\widetilde P}_b/ d{\widetilde P}_0)|{\mathcal C}_{A_n}^Y= {\widetilde
L}_{A_n}(b)= {\widetilde L}_{S_n}(b) V_n$,
with $
\log{V_n}=\int_{S_n}^{A_n} f(Y_u) \,dY_u - \int_{S_n}^{A_n} \frac
{1}{2} f^2(Y_u) \,du$.
Hence,\vspace*{2pt}
$
d{\widetilde P}_b|{\mathcal C}_{A_n}^Y/dN({\widetilde P}_b|{\mathcal
C}_{S_n}^Y)= V_n$.
By the Pinsker inequality (\hyperref[secappen]{Appendix}) and Lemma~\ref
{rhorhobar}, we have
\begin{eqnarray*}
\bigl\|N\bigl({\widetilde P}_b|{\mathcal C}_{S_n}^Y
\bigr)- {\widetilde P}_b|{\mathcal C}_{A_n}^Y
\bigr\|_{\mathrm{TV}}&=& \frac{1}{2} \int_{\Omega} d{\widetilde
P}_0\bigl|{\widetilde L}_{S_n}(b)-{\widetilde
L}_{A_n}(b)\bigr|
\\
&\le& \sqrt{K\bigl({\widetilde P}_b|{
\mathcal C}_{A_n}^Y,N\bigl({\widetilde P}_b|{
\mathcal C}_{S_n}^Y\bigr)\bigr)/2},
\\
K\bigl({\widetilde P}_b|{\mathcal C}_{\tau_n},N({\widetilde
P}_b|{\mathcal C}_{S_n})\bigr) &=& E_{{\widetilde P}_b|{\mathcal C}_{A_n}^Y} \int
_{S_n}^{A_n} \frac{1}{2} f^2(X_u)
\,du
\\
&\le& \frac{K^2}{2 \sigma_0^4} E_{{\widetilde P}_b} |A_n - {
\widebar{A}_n}| \le\frac{K^2}{ \sigma_0^4}c nh_n^2.
\end{eqnarray*}
Using that
%We obtain the first inequality using that
$
\delta({\widetilde{\mathcal E}}_0^{S_n},{\widetilde{\mathcal E}}_0^{A_n})\le
\sup_{b\in{\mathcal F}_K}\|N({\widetilde P}_b|{\mathcal C}_{S_n}^Y)-
{\widetilde P}_b|{\mathcal C}_{A_n}^Y\|_{\mathrm{TV}}
$ yields the first inequality.
We proceed analogously for the other one.

\textit{Proof of inequality}~(\ref{DistanceTV}).
These experiments have the same sample space and are, respectively,
associated with the distributions ${\widetilde P}_{b}$ (resp., $\widetilde{
Q}_{b}$)
on $C({\mathbb R}^+, {\mathbb R})$ of $(\zeta_u, u \ge0)$ given by
(\ref{defY}) [resp., $(\bar\zeta_u, u \ge0)$ given by~(\ref{Ybar})].
%The distribution ${\widetilde P}_{b}$ (resp. ${\widetilde Q}_{b}$) is
%associated with the diffusion process with drift $f$ (resp. the
%diffusion type process with drift ${\bar f}(v,y)$) with the same
%initial condition.
Hence,
\[
\Delta\bigl({\widetilde{\mathcal E}}_{0}^{{\widebar A}_n}, {\widetilde{
\mathcal G}}_{0}^{{\widebar A}_n}\bigr) \le \sup_{b \in{\mathcal F}}
\| {\widetilde P}_{b}/_{{\mathcal C}_{{\widebar
A}_n^Y}} -\widetilde{
Q}_{b}/_{{\mathcal C}_{{\widebar A}_n^Y}} \|_{\mathrm{TV}}=\Delta _0
\bigl({\widetilde{\mathcal E}}_{0}^{{\widebar A}_n}, {\widetilde{\mathcal
G}}_{0}^{{\widebar A}_n}\bigr). %
\]
Using the bound of Proposition~\ref{Pinsk} yields
\begin{eqnarray*}
2\| {\widetilde P}_{b}/_{{\mathcal C}_{{\widebar A}_n}^Y} - {\widetilde
Q}_{b}/_{{\mathcal C}_{{\widebar A}_n}^Y}\|_{\mathrm{TV}}^2 &\le& K({
\widetilde P}_{b}/_{{\mathcal C}_{{\widebar A}_n}^Y}, {\widetilde Q}_{b}/_{{\mathcal C}_{{\widebar A}_n}^Y})
\\
&=& E_{{\widetilde P}_b/_{{\mathcal C}_{{\widebar A}_n}^Y} } \biggl( \int_0^{{\widebar A}_n} \bigl(
f(Y_v)-{\bar f}(v,Y)\bigr) ^2 \,dv \biggr). %
\end{eqnarray*}
Setting $T_{i}=T_i(Y)$ and using that, for $i=1, \ldots, n$,
$T_{i}=T_i(Y)={\widebar A}_{t_i}(Y)$ [see~(\ref{lemTiYTi})] and that
$f$ is Lipschitz with constant $L$ [see~(\ref{f})], we get
\begin{eqnarray*}
\int_0^{{\widebar A}_n} \bigl( f(Y_v)-{\bar
f}(v,Y)\bigr) ^2 \,dv &=& \sum_{i=1}^n
\int_{ T_{i-1}}^{ T_{i}} \bigl(f(Y_v) -
f(Y_{T_{i-1}} )\bigr)^2 \,dv
\\
&\le& {L^2}\sum
_{i=1}^n \int_{ T_{i-1}}^{ T_{i}}
(Y_v- Y_{T_{i-1}} )^2 \,dv. %
\end{eqnarray*}
Under ${\widetilde P}_b$, $Y_v- Y_{T_{i-1}}= \int_{T_{i-1}}^v f(Y_u) \,du +
B_v -B_{T_{i-1}}$,
with $(B_v)$ Brownian motion.~So
\[
(Y_v- Y_{T_{i-1}}) ^2\le2  \biggl[  \biggl(\int_{T_{i-1}}^v f(Y_u) \,du \biggr)^2 + (B_v -B_{T_{i-1}})^2 \biggr].
\]
This yields $
\int_0^{{\widebar A}_n} ( f(Y_v)-{\bar f}(v,Y)) ^2 \,dv \le2 L^2 (R_1+R_2)$,
with
\[
R_1= \sum_{i= 1}^n \int_{ T_{i-1}}^{ T_{i}}  \biggl(\int_{T_{i-1}}^u
f(Y_v) \,dv  \biggr)^2 \,du,\qquad
 R_2= \sum_{i =1}^n \int_{ T_{i-1}}^{
T_{i}}  (B_u - B_{T_{i-1}}  )^2 \,du.
\]
Using~(\ref{f}) and $T_i-T_{i-1} \leq\sigma_1^2h_n$ by~(\ref{defTi}),
%by~(\ref{defTi}),
\[
R_1\le\frac{K^2}{\sigma_0^4} \sum_{i= 1}^n (T_{i}- T_{i-1})^3 \le
\frac{K^2}{\sigma_0^4} n\bigl(\sigma_1^2 h_n\bigr)^3.
\]
For the second term, using definition~(\ref{bbari}),
\begin{eqnarray*}
E_{{\widetilde P}_b}(R_2)&=& E_{{\widetilde P}_b} \Biggl(\sum
_{i =1}^n \int_{
T_{i-1}}^{ T_{i}}
(B_u - B_{T_{i-1}} )^2 \,du\Biggr) \le\sum
_{i=1}^{n} \int_0^{\sigma_1^2 h_n}
E_{{\widetilde P}_b}\bigl({\widebar B}_v^{(i)}
\bigr)^2 \,dv
\\
&=& n \frac{(\sigma_1^2h_n)^2}{2}. %
\end{eqnarray*}
Thus, the result follows from
\[
K({\widetilde P}_{b}/_{{\mathcal C}_{{\widebar A}_n^Y}}, {\widetilde
Q}_{b}/_{{\mathcal C}_{{\widebar A}_n^Y}}) \le2 L^2 \biggl(
\frac
{K^2}{3\sigma_0^4} n\bigl(\sigma_1^2 h_n
\bigr)^3+ n \frac{(\sigma
_1^2h_n)^2}{2} \biggr).
\]%upqed
Joining~(\ref{propsuffi}),~(\ref{DistanceTV}) and~(\ref{lastineq})
completes the proof of Proposition~\ref{proprandom}.
\end{pf*}

\begin{appendix}\label{sec6} \label{secappen}
%s6 #&#
\section*{Appendix}
Let us recall properties of the Le Cam deficiency distance $\Delta$.
%All measurable spaces are supposed to be Polish metric spaces equipped
%with their Borel $\sigma$-algebras.
Consider two statistical experiments ${\mathcal E}=(\Omega, {\mathcal
A}, (P_f)_{f \in{\mathcal F}})$ and
${\mathcal G}=({\mathcal X}, {\mathcal C}, (Q_f)_{f \in{\mathcal F}})$
and assume that the families $(P_f)_{f \in{\mathcal F}}$, $(Q_f)_{f
\in{\mathcal F}}$ are dominated. A Markov kernel $M(\omega, dx)$ from
$(\Omega, {\mathcal A})$ to $({\mathcal X}, {\mathcal C})$ is a
mapping from $\Omega$ into the set of probability measures on
$({\mathcal X}, {\mathcal C})$ such that, for all $C \in{\mathcal C}$,
$\omega\rightarrow M(\omega,C)$ is measurable on $(\Omega, {\mathcal
A})$, and for all $\omega\in\Omega$, $M(\omega, dx)$ is probability
measure on $({\mathcal X}, {\mathcal C})$. The image $MP_f$ of $P_f$
under $M$ is defined by
$
MP_f(C)= \int_{\Omega} M(\omega, C) \,dP_f(\omega)$.
The experiment $M{\mathcal E}=({\mathcal X}, {\mathcal C}, (MP_f)_{f
\in{\mathcal F}})$ is called a randomization of ${\mathcal E}$ by the
kernel $M$. If the kernel is deterministic, that is, for $T\dvtx (\Omega,{\mathcal A}) \rightarrow({\mathcal X}, {\mathcal C})$ a random
variable, $T(\omega,C)=1_{C}(T(\omega))$, the experiment $T{\mathcal
E}$ is called the image experiment by $T$.
%
%de6.1 #&#
\begin{defi}
 $\Delta({\mathcal E},{\mathcal G})= \max\{\delta
({\mathcal E},{\mathcal G}),\delta({\mathcal G},{\mathcal E})\}$ where
$\delta({\mathcal E},{\mathcal G})=\break  \inf_{M \in{\mathcal M}_{\Omega\dvtx {\mathcal X}}} \sup_{f \in{\mathcal F}}\|MP_f- Q_f\|_{\mathrm{TV}}$,
$\|\cdot\|_{\mathrm{TV}}$ is the total variation distance and ${\mathcal M}_{\Omega\dvtx {\mathcal X}}$ the set of Markov kernels from $(\Omega, {\mathcal
A})$ to $({\mathcal X}, {\mathcal C})$.
\end{defi}

When $\Delta({\mathcal E},{\mathcal G})=0$, the two experiments are
said to be equivalent. When the experiments have the same sample space:
$(\Omega,{\mathcal A})=({\mathcal X}, {\mathcal C})$, it is possible
to define
$
\Delta_0({\mathcal E},{\mathcal G})= \sup_{f \in{\mathcal F}}\|P_f-
Q_f\|_{\mathrm{TV}}$, which satisfies
$
\Delta({\mathcal E},{\mathcal G}) \le\Delta_0({\mathcal E},{\mathcal G})$.
Consider an asymptotic framework $\varepsilon\rightarrow0$ and
families of experiments ${\mathcal E}^{\varepsilon}=(\Omega
^{\varepsilon},
{\mathcal A}^{\varepsilon}, (P_f^{\varepsilon})_{f \in{\mathcal F}})$,
${\mathcal G}^{\varepsilon}=({\mathcal X}^{\varepsilon}, {\mathcal
C}^{\varepsilon}, (Q_f^{\varepsilon})_{f \in{\mathcal F}})$,
${\mathcal B}^{\varepsilon}\subset{\mathcal A}^{\varepsilon}$ a
$\sigma$-algebra.
%, $T^{\varepsilon}\dvtx (\Omega^{\varepsilon},{\mathcal A}^{\varepsilon})
%a statistic.
%
%de6.2 #&#
\begin{defi} The families ${\mathcal E}^{\varepsilon}$, ${\mathcal
G}^{\varepsilon}$ are asymptotically equivalent as $\varepsilon
\rightarrow0$ if $\Delta({\mathcal E}^{\varepsilon},{\mathcal
G}^{\varepsilon})$ tends to $0$. The $\sigma$-algebra ${\mathcal
B}^{\varepsilon}$ is asymptotically sufficient if
$\Delta({\mathcal E}^{\varepsilon},{\mathcal E}^{\varepsilon
}/_{{\mathcal B}^{\varepsilon}} )$ tends to $0$, where ${\mathcal
E}^{\varepsilon}/_{{\mathcal B}^{\varepsilon}}$ is the restriction of
${\mathcal E}^{\varepsilon}$ to ${\mathcal B}^{\varepsilon}$.\vadjust{\goodbreak}
\end{defi}
We state now two auxiliary results used in proofs.
%
%pr6.1 #&#
\begin{prop}\label{propauxi}
Let $(B_t, t\ge0)$ be a Brownian motion
with respect to a filtration $({\mathcal F}_t, t\ge0)$ (satisfying the
usual conditions) and let $\tau$ be a positive ${\mathcal
F}_0$-measurable random variable. Then $(W(t)= \frac{1}{\sqrt{\tau
}}B_{\tau t}, t\ge0)$ is a standard Brownian motion, independent of
${\mathcal F}_0$.
\end{prop}
This result follows from a straightforward application of Paul L\'evy's
characterisation of the Brownian motion [see, e.g., \citet{KA}].
%0$, $\tau t$ is a $({\mathcal F}_t)$- stopping time.
%By the optional sampling theorem, we deduce that the processes $(B_{
%$(B^{2}_{\tau t} - \tau t, t\ge0)$ are local martingales with respect
%to the filtration
%$({\mathcal F}_{\tau t})$, with continuous sample paths, null at $0$.
%As $\tau$ is ${\mathcal F}_0$-measurable, the same holds for the two
%processes $(W(t))$ and $(W^2(t) -t)$. Thus, by Paul L\'evy's
%characterization, we deduce that $(W(t))$ is a Brownian motion, with
%respect to the same filtration. Thus, $(W(t))$ is independent of ${
Next, we recall the first Pinsker inequality [see, e.g., \citet{Tsyb}] for the total variation distance between probability measures.
Let $({\mathcal X}, {\mathcal A})$ be a measurable space, $P, Q$ two
probability measures on $({\mathcal X}, {\mathcal A})$, $\nu$ a
$\sigma$-finite measure on $({\mathcal X}, {\mathcal A})$ such that $P
\ll\nu$, $Q \ll\nu$ and set $p=dP/d\nu, q=dQ/d\nu$. The total
variation distance between $P$ and $Q$ is defined by:
$
\|P-Q\|_{\mathrm{TV}}= \sup_{A \in{\mathcal A}}|P(A)-Q(A)|= \frac{1}{2}\int
|p-q|\,d\nu$.
The Kullback divergence of $P$ w.r.t. $Q$ is $K(P,Q)= \int\log{\frac
{dP}{dQ}} \,dP$ if $P\ll Q$, $=+\infty$ otherwise.
%
%pr6.2 #&#
\begin{prop}\label{Pinsk}
$
\|P-Q\|_{\mathrm{TV}} \le\sqrt{K(P,Q)/2}$.
\end{prop}
The remarkable feature of this inequality is that the left-hand side is
a symmetric quantity whereas the right-hand side is not. The noteworthy
consequence is that it is possible to choose, for the right-hand side,
$K(P,Q)$ or $K(Q,P)$. The Pinsker inequality is particularly useful
when $P,Q$ are associated with diffusion type processes. Let $P$
(resp., $Q$) be the distribution $C({\mathbb R}^+, {\mathbb R})$ of the
diffusion type process $d\xi_t=p(t, \xi_.) \,dt + dW_t$ with
predictable drift $p(t, X_.)$ [resp., $d\eta_t=q(t, \eta_.)\,dt + dW_t$
with drift $q(t,X_.)$] and constant diffusion coefficient equal to $1$,
with the same initial condition $\xi_0=\eta_0$. Let $T=T(X_.)$ be a
finite stopping time under $P$ and $Q$. Then the Girsanov formula
stopped at $T$ yields
[with $(X_v)$ the canonical process of $C({\mathbb R}^+, {\mathbb R})$]
\[
\frac{d P_T}{dQ_T}=\exp{ \biggl(\int_0^{T}
\bigl(p(s,X_.)-q(s,X.) \bigr) \,dX_s- \frac{1}{2} \int
_0^{T} \bigl(p^2(s,X_.)-q^2(s,
X_.) \bigr)\,ds \biggr)},
\]
where $P_T= P/_{{\mathcal C}_{T} }, Q_T= Q/_{{\mathcal C}_{T} }$ are
the restriction of $P, Q$ to the $\sigma$-field ${\mathcal C}_{T} $.
Hence, using that under $P\, dX_t-p(t,X_.)\,dt = dB_t $, with $(B_t)$ a
Brownian motion, yields
$
K(P_T, Q_T)
=(1/2) E_{P }( \int_0^T (p(s,X_.)-q(s,X_.) ) ^2 \,ds )$.
\end{appendix}


%suskaldyti doi

% imsref loaded by linak, 2014-04-16 12:31:52

\printaddresses

\end{document}